\documentclass[12pt]{amsart}
\usepackage{a4wide,amssymb,mathrsfs}
\usepackage[2cell,arrow,matrix,curve]{xy}

\let\x\times

\let\Ph\varPhi

\let\ka\varkappa
\let\d\partial

\let\ge\geqslant
\let\normal\triangleleft

\def\sJ{{\mathscr J}}
\def\1{^{-1}}

\def\xto#1{\xrightarrow[]{#1}}

\DeclareMathOperator\cro{cr} \DeclareMathOperator\id{id}

\DeclareMathOperator\Sym{S}

\def\set#1{\left\{#1\right\}}
\def\brk#1{\langle#1\rangle}

\newtheorem{Pro}{Proposition}%[section]
\newtheorem{Le}[Pro]{Lemma}
\newtheorem{The}[Pro]{Theorem}
\newtheorem{Co}[Pro]{Corollary}
\theoremstyle{definition}
\newtheorem{De}[Pro]{Definition}
\theoremstyle{remark}

\newtheorem{Rem}[Pro]{Remark}

\def \Im{\mathop{\sf Im}\nolimits}
\def \Ker{\mathop{\sf Ker}\nolimits}
\def \Hom{\mathop{\sf Hom}\nolimits}

\newcommand{\QZ}{{\mathbb Q}/{\mathbb Z}}

\def\C{{\cal C}}

\def\Nil{{\sf Nil}}
\def\lin{{\sf Lin}}
\def\qgr{{\sf Quad}}
\def\qgrs{{\sf Quad}^{\Sigma}}
\def\nilp{^{\mathrm{nil}}}
\def\Sets{{\sf Sets}}
\let\x\times

\let\tp\otimes
\def\La{\Lambda}

\let\t\otimes

\def\Cok{\mathop{\sf Coker}\nolimits}

\let\cal\mathscr

\let\ge\geqslant
\let\le\leqslant

\let\x\times
\let\ox\otimes

\let\into\rightarrowtail
\let\onto\twoheadrightarrow

\let\Ph\varPhi

\def\xto#1{\xrightarrow[]{#1}}

\def\ab{^{\mathrm{ab}}}

\def\brk#1{\left\langle#1\right\rangle}
\def\set#1{\left\{#1\right\}}

\def\id{{\mathrm{Id}}}

\def\op{^{\mathrm{op}}}
\def\1{^{-1}}

\def\Znil{{{\mathbb Z}_{\mathrm{nil}}}}

\def\e{_{\mathrm e}}
\def\ee{_{\mathrm{ee}}}

\def\LL{{\mathbb L}}

\def\ZZ{{\mathbb Z}}

\def\ssq{{\mathsf{Simpl(SG)}}}
\def\Square{{\mathsf{SG}}}

\def\ld{_{\mathrm l}}
\def\rd{_{\mathrm r}}

\def\cref#1#2#3{\left(#2\right.\left|\ #3\right)_{#1}}

\def\centr{{\mathscr Z}}

\def\PSG{{\sf PSG}}
\def\Cos{{\sf Cos}}
\def\SG{{\sf SG}}
\def\Sets{{\sf Sets}}
\def\QZ{{\sf QZ}}

\def\Ab{{\mathsf{Ab}}}
\def\Gr{{\sf Groups}}

\def\Id{{\sf Id}}

\def\tl{{\underline\circledcirc}}
\let\tr\circledcirc
\def\tb{{\bar\otimes}}
\let\un\upsilon
\def\ttl{\mathop{{\sf Tor}^{\odot}}\nolimits}
\begin{document}

\title{Quadratic algebra of square groups}

\author
{H.-J. Baues}

\address{Max-Planck-Institut f\"ur Mathematik, 53111 Bonn, Germany}

\author{M. Jibladze}

\address{Razmadze Mathematical Institute, Tbilisi 0193, Georgia}

\author{T. Pirashvili}

\address{Razmadze Mathematical Institute, Tbilisi 0193, Georgia}

\thanks{The second and third authors are grateful to the Max-Planck-Insitut f\"ur Mathematik, Bonn, where this work was written, for hospitality}

\begin{abstract} Square groups are quadratic analogues of abelian groups. Many properties of abelian groups are shown to
hold for square groups. In particular, there is a symmetric monoidal tensor product of square groups generalizing the
classical tensor product.
\end{abstract}
%\hfill(Version of \today)

\maketitle
%\begin{keyword}{Quadratic functors, Square groups, symmetric monoidal category}
%\MSC 18D10, 18G50, 55Q99.
%\end{keyword}
%\maketitle
%\part{Description of main results}

%\section{Introduction}
There is a long-standing problem of algebra to extend the symmetric monoidal
structure of abelian groups, given by the tensor product, to a non abelian
setting, see for example \cite{BL}. In this paper we show the somewhat
surprising fact that such an extension is possible. Morover our non abelian
tensor product remains even right exact and balanced. We describe the new
non-abelian tensor product in the context of quadratic algebra which extends
linear algebra.

``Linear algebra'' is the algebra of rings and modules. A ring is a monoid in
the symmetric monoidal category of abelian groups
$$
(\Ab,\ox,\ZZ).
$$
The monoidal structure is given by the tensor product
of abelian groups, with the group of integers $\ZZ$ as the unit object.
Moreover a module is an object in $\Ab$ together with an action of such a monoid.

In ``quadratic algebra'' abelian groups are replaced by square groups. In
fact, if one considers endofunctors of the category of groups which preserve
filtered colimits and reflexive coequalizers, then abelian groups can be
identified with linear endofunctors and square groups can be identified with
quadratic endofunctors, \cite{square}. The abelian group $\ZZ$ corresponds to
the linear endofunctor which carries a group $G$ to its abelization
$G\ab=G/[G,G]$. The square group $\Znil$ corresponds to the quadratic
endofunctor which carries $G$ to the class two nilpotent group
$G\nilp=G/[G,[G,G]]$. The category $\Square$ of square groups
contains the category $\Ab$ of abelian groups as a full subcategory since a
linear endofunctor is also quadratic. Therefore the question arises whether
the symmetric monoidal structure of $\Ab$ extends to a symmetric monoidal
structure of $\Square$. The main purpose of this paper is the proof that this
is, in fact, the case.

Let $G$ and $H$ be
(additively written) groups. One can consider the group with generators $g*h$
for $g\in G$ and $h\in H$, subject to the following relations:
\begin{align*}
g*(h_1+h_2)&=g*h_1+g*h_2,\\
(g_1+g_2)*h&=g_1*h+g_2*h.
\end{align*}
It is well known and easy to prove that this group is isomorphic to $G\ab\t
H\ab$. Thus the naive definition of the tensor product of
nonabelian groups does not really provide a new object. A more sophisticated
tensor product was constructed by Brown and Loday \cite{BL}. However, their
tensor product does not define a  symmetric monoidal structure.

The definition of the tensor product of square groups relies on
the structure
$(H,P)$ of a square group $M$. In $M\tl N$ one has among others the relations
\begin{align*}
x\tl(y_1+y_2)&=x\tl y_1+x\tl y_2,\\
(x_1+x_2)\tl y&=x_1\tl y+x_2\tl y+\cref H{x_2}{x_1}\tb H(y)
\end{align*}
which replace the naive relations above. It is a somewhat surprising fact that
these relations for $M\tl N$ lead to a symmetric monoidal structure of the
category of square groups extending the tensor product of abelian groups.
We show:

{\bf Theorem}. \emph{
There is a tensor product of square groups $M$, $N$ denoted by $M\tl N$ such
that
$$
(\Square,\tl,\Znil)
$$
is a symmetric monoidal category. Moreover, if $M$ and $N$ are abelian groups
then
$$
M\tl N=M\ox N.
$$
}
%\end{The}

%We give three equivalent constructions $M\tr N\cong M\tl N\cong M\odot N$.
A
monoid $R$ in $(\Square,\tl,\Znil)$ is termed a quadratic ring. An $R$-quadratic module is a square group with an action of the monoid $R$. This leads to the
wide area of quadratic algebra generalizing classical linear algebra. For example
we describe in this paper the Tor-exact sequence for square groups. Also we study
various special classes of square groups, like abelian square groups, quadratic
$\ZZ$-modules and free square groups.

The category $\Square$ has another (very non-symmetric) monoidal category
structure $\square$ induced by composition of endofunctors (see \cite{square}
or Section \ref{kvfunqtorebi}). In Theorem \ref{jibladzestrange} we describe,
by means of abelian groups with cosymmetry, a subcategory of $\Square$
on which the products $\square$ and $\tl$ coincide.

As we will see in sequel publications quadratic rings and related quadratic
pair algebras play an important r\^ole in secondary homotopy theory as well
as in the theory of Mac Lane cohomology of rings (see \cite{JP}, \cite{PW},
\cite{shukla} and Chapter 13 of \cite{HC}). Namely, quadratic pair algebras
are natural objects representing classes in the third dimensional Mac Lane
cohomology. Moreover the secondary homotopy groups of each ring spectrum form
a quadratic pair algebra \cite{Bm}. In particular, the sphere spectrum yields
a quadratic pair algebra encoding all its secondary homotopy structure like
triple Toda brackets. Also in order to study these examples it is necessary
to develop the quadratic algebra of square groups.

Concerning (symmetric) monoidal categories and (symmetric) monoidal functors
we use the terminology following \cite{js}. In particular, a lax monoidal
functor is a functor $F$ together with coherent morphisms $\phi_{A,B}:F(A)\t
F(B)\to F(A\t B)$. Moreover $F$ is a monoidal functor if $\phi_{A,B}$ are
isomorphisms for all $A$ and $B$.

\section{The monoidal category of square groups}\label{11}
In this section we recall the notion of square group (see also Section
\ref{kvfunqtorebi} below) and we give an explicit construction of the tensor
product of square groups. We formulate our main results concerning symmetric
monoidal structure and right exactness of this tensor product.

Let $G$ be an additively written group and $A$ be an abelian group. We call a
map $f:G\to A$ \emph{quadratic} if for any $x,y\in G$ the \emph{cross-effect}
$$
\cref fxy:=f(x+y)-f(x)-f(y)
$$
is linear in $x$ and $y$, that is $\cref f-y$ and $\cref fx-$ are homomorphisms
$G\to A$.

\begin{De}\label{sgdef} A \emph{square group} is a diagram
$$
M = (\xymatrix{M\e \ar[r]^H & M\ee \ar[r]^P & M\e})
$$
where the $\mathrm{ee}$-level $M\ee$ is an abelian group and the
$\mathrm{e}$-level $M\e$ is a group. Both groups are written additively.
Moreover $P$ is a homomorphism and $H$ is a quadratic map. In addition the
following identities
\begin{align*}
\cref H{Pa}y&=0=\cref Hx{Pb},\\
P\cref Hxy&=-x-y+x+y,\\
PHP(a)&=2P(a)
\end{align*}
 are satisfied for all $x,y\in M\e$ and $a,b \in M\ee$.
\end{De}
As an example we have the square group
$$
\Znil=(\xymatrix{\ZZ \ar[r]^H& \ZZ\ar[r]^P&\ZZ})
$$
with $H(n)=\binom n2=\frac{n(n-1)}{2}$ and $P=0$. Let $\Square $ be the
category of square groups.

In any square group $M$ the image of $P$ is a normal subgroup containing the
commutator subgroup (see \cite{square}, or Section \ref{squag}), thus
$\Cok(P)$ is a well-defined abelian group. Let $\bar{x}\in \Cok(P)$ be the
element represented by $x\in M\e$. The cross effect of $H$ induces a
homomorphism (see \cite{square} or Corollary \ref{29} below)
$$\cref H-- :\Cok(P)\t \Cok(P)\to M\ee.$$
Moreover, there is a well-defined
homomorphism
$$
\Delta:\Cok(P)\to M\ee
$$
given by $\Delta(\bar{x})=HPH(x)+H(x+x)-4H(x)$ (see \cite{square}, or Corollary
\ref{29}). Furthermore, the map
$$
T=HP-\Id:M\ee\to M\ee
$$
is an endomorphism of the abelian group $M\ee$ with $T^2=\Id$ (see
\cite{square}, or Proposition \ref{sqidentities}). Sometimes we write
$P^M,H^M,\Delta^M,T^M$ in order to make clear the r\^ole of $M$.

The aim of this work is to introduce a symmetric monoidal category structure
$\tl$ on the category $\Square$ of square groups.

\begin{De}\label{tldef}
For two square groups $M$, $N$ we introduce the tensor product $M\tl N$ which
is a square group defined as follows. The group $(M\tl N)\e$ is given by
generators of the form $x\tl y$ for $x\in M\e$, $y\in N\e$, and $a\tb b$ for
$a\in M\ee$, $b\in N\ee$, subject to the relations
\begin{enumerate}
\item\label{bice0} the symbol $a\tb b$ is bilinear and central in $(M\tl N)\e$;
\item\label{ledi0} $x\tl(y_1+y_2)=x\tl y_1+x\tl y_2$;
\item\label{ridi0} $(x_1+x_2)\tl y=x_1\tl y+x_2\tl y+\cref H{x_2}{x_1}\tb H(y)$;
\item\label{rip0} $x\tl P(b)=\cref Hxx\tb b$;
\item\label{lep0} $P(a)\tl y=a\tb\Delta(y)$;
\item\label{TT0} $T(a)\tb T(b)=-a\tb b$.
\end{enumerate}
Next the abelian group $(M\tl N)\ee$ is defined to be $M\ee\ox N\ee$. The
homomorphism
$$
P:(M\tl N)\ee\to(M\tl N)\e
$$
is given by
$$
P(a\ox b)=a\tb b
$$
In order to define $H$, we first observe that
$$
\Cok(P^{M\tl N})=\Cok(P^M)\t \Cok(P^N).
$$
Therefore we have the following homomorphism
\begin{multline*}
\rho:\Cok(P^{M\tl N}) \t \Cok(P^{M\tl N})\\
=\Cok(P^M)\t \Cok(P^N)\t\Cok(P^M)\t\Cok(P^N)\to M\ee\t N\ee
\end{multline*}
with $\rho(\bar a\t \bar b\t \bar a'\t \bar b')=\cref Ha{a'}\t\cref Hb{b'}$. Now
$$
H:(M\tl N)\e\to(M\tl N)\ee
$$
is the unique quadratic map with the map $\rho$ as its cross-effect
satisfying
$$
H(x\tl y)=\cref Hxx\ox H(y)+H(x)\ox\Delta(y)
$$
and
$$
H(a\tb b)=a\ox b-T(a)\ox T(b).
$$
\end{De}

The following is the main result of the paper:
\begin{The}\label{tlmonoidal}
The tensor product of square groups gives rise to a well-defined bifunctor
$$
-\tl-:\Square\x\Square\to\Square
$$
which equips the category $\Square$ with a symmetric monoidal structure,
with the unit object $\Znil$. The associativity and commutativity isomorphisms
on the $\mathrm{ee}$-level are the usual isomorphisms for the tensor product of
abelian groups, while on the $\mathrm{e}$-level the isomorphism
$((M\tl N)\tl K)\e\cong
(M\tl (N\tl K))\e$ is given by
\begin{align*}
(x\tl y)\tl z&\mapsto x\tl (y\tl z)\\
(a\tb b)\tl z&\mapsto a\tb (b\t \Delta(z))\\
(a\t b)\tb c&\mapsto a\tb (b\t c)
\end{align*}
and the isomorphism $(M\tl N)\e\cong (N\tl M)\e$ is given by
\begin{align*}
x\tl y&\mapsto y\tl x-H(y)\tb TH(x)\\
a\tb b&\mapsto b\tb a.
\end{align*}
 Under the associativity
isomorphism the element $x\tl (b\tb c)$ corresponds to $(x|x)_H\tb (b\t c)$.
\end{The}

Proof of this result occupies Section \ref{tl901}.

The category of square groups with $M\ee=0$ is equivalent to the category of
abelian groups. Thus we identify the category $\Ab$ with this subcategory of
$\Square$. Then the restriction of $\tl$ to $\Ab$ coincides with the usual
tensor product of abelian groups. The following generalizes some well-known
properties of the tensor product of abelian groups.

\begin{Pro}\label{tlright}

For any square group $A$ the tensor product $A\tl - :\Square\to \Square$
preserves filtered colimits, reflexive coequalizers and finite
products. It is right exact and balanced, that is
for any short exact sequence of square groups
$$
0\to B_1\xto{\mu} B\xto{\sigma} B_2\to 0
$$
the induced sequence
$$
A\tl B_1\to A\tl B\to A\tl B_2\to 0
$$
is exact and the first arrow $A\tl \mu$ is a monomorphism provided
$A$ is a projective object in the category $\Square$.
\end{Pro}

The proof of this result is given in Section \ref{rexacttl}.

The functor
$A\tl - :\Square\to \Square$ does not preserves coproducts. However the following
result is true. For an abelian group $A$ we define the square group $A^\t$ by
$$(A^\t)\e=A, \ \ (A^\t)\ee=A\oplus A,$$
where $P(a,b)=a+b$ and $H(a)=(a,a)$.
\begin{Pro}\label{tlco} Let $A,B,M$ be square groups. Then one has the short exact
sequence of square groups
$$
0\to(M\ee\t\Cok(P^A)\t\Cok(P^B))^\t\to M\tl(A\vee B)\to(M\tl A)\x(M\tl
B)\to0.
$$
Here $\vee$ denotes the coproduct in the category of square groups.
\end{Pro}
The proof of this result is given at the end of Section \ref{tortl}.

\section{Symmetric definition of the tensor product}
The tensor product $M\tl N$ in Definition \ref{tldef} is the \emph{right linear
version}. There is also a \emph{left linear version} $M\tr N$ defined below.
Moreover we introduce a \emph{symmetric version} $M\odot N$ and we show that
there are natural isomorphisms
$$M\tr N\cong M\odot N\cong M\tl N.$$
The different versions of the tensor product are identified in this way. The
symmetry of the tensor product is most apparent in $M\odot N$ where, however,
redundant generators are needed. The non-symmetric versions $M\tl N$ and
$M\tr N$ have the advantage of a smaller set of generators. Most calculations
in the paper are using the right linear version $M\tl N$.

\begin{De}\label{trdef}
For square groups $M$, $N$ let $M\tr N$ be the
 square group with $(M\tr N)\ee=M\ee \t N\ee$ and $(M\tr N)\e$ given
by generators $x\tr y$ for $x\in M\e$, $y\in N\e$, and $a\tb b$ for $a\in
M\ee$, $b\in N\ee$, subject to the relations
\begin{enumerate}
\item\label{rbice0} the symbol $a\tb b$ is bilinear and central in $(M\tr N)\e$;
\item\label{rledi0} $x\tr(y_1+y_2)=x\tr y_1+x\tr y_2+H(x)\tb \cref H{y_2}{y_1}$;
\item\label{rridi0} $(x_1+x_2)\tr y=x_1\tr y+x_2\tr y$;
\item\label{rrip0} $x\tr P(b)=\Delta(x)\tb b $;
\item\label{rlep0} $P(a)\tr y=a\tb\cref Hyy$;
\item\label{rTT0} $T(a)\tb T(b)=-a\tb b$.
\end{enumerate}
Here the homomorphism $P$ is given as in Definition \ref{tldef} and $H$ is
the quadratic map with the cross effect $\rho$
as in Definition \ref{tldef} and
$$
H(x\tr y)=\Delta(x)\t H(y)+ H(x)\t
\cref Hyy,
$$
$$
H(a\tb b)=a\ox b-T(a)\ox T(b).
$$
\end{De}
Next we introduce the symmetric version of the tensor product.

\begin{De}\label{symdef}
For two square groups $A$, $B$ we define their tensor product $A\odot B$ which
is again a square group defined as follows. The group $(A\odot B)\e$ is defined
by generators of the form $x\tl y$, $x\tr y$ for $x\in A\e$, $y\in B\e$ and $a\tb b$ for
$a\in A\ee$, $b\in B\ee$, subject to the relations
\begin{enumerate}
\item\label{sbice} the symbol $a\tb b$ is bilinear and central in $(A\tl B)\e$;
\item\label{sledi} $x\tl(y_1+y_2)=x\tl y_1+x\tl y_2$;
\item\label{sridi} $(x_1+x_2)\tr y=x_1\tr y+x_2\tr y$;
\item\label{slep} $P(a)\tl y=a\tb\Delta(y)$;
\item\label{srip} $x\tr P(b)=\Delta(x)\tb b$;
\item\label{sTT} $T(a)\tb T(b)=-a\tb b$;
\item\label{ssym} $x\tl y-x\tr y=H(x)\tb TH(y)$.
\end{enumerate}
Next the abelian group $(A\odot B)\ee$ is defined to be $A\ee\ox B\ee$. The
homomorphism
$$
P:(A\odot B)\ee\to(A\odot B)\e
$$
is given by
$$
P(a\ox b)=a\tb b
$$
and the map
$$
H:(A\odot B)\e\to(A\odot B)\ee
$$
is the unique quadratic map with
\begin{align*}
H(x\tl y)&=(x|x)_H\ox H(y)+H(x)\ox\Delta(y),\\
H(x\tr y)&=\Delta(x)\ox H(y)+H(x)\ox(y|y)_H
\end{align*}
and
$$
H(a\tb b)=a\ox b-T(a)\ox T(b)
$$
such that its cross-effect coincides with the bilinear map $\rho$ in Definition
\ref{tldef}.
\end{De}

\begin{Pro}\label{tl=tr=to}
For any square groups $A$, $B$ the above data define square groups $A\odot B$
and $A\tr B$ and both of them are isomorphic to $A\tl B$.
\end{Pro}

The proof of this fact is given in Section \ref{totr}.

\begin{Co}
The symmetry isomorphism
$$
\tau(A,B):A\odot B\to B\odot A,
$$
corresponding to the symmetry for $A\tl B$ under the above isomorphism, is given by
\begin{align*}
x\tl y&\mapsto y\tr x,\\
x\tr y&\mapsto y\tl x.\\
a\tb b &\mapsto b\tb a.
\end{align*}
\end{Co}

\begin{Rem}
The notation above is chosen to be compatible with the notation for
\emph{exterior cup products} $f\#g$ and $f\underline\#g$ in topology, see
\cite{comca}, \cite{dreck}. Here $f\#g$, being left linear, corresponds to
$x\tr y$ and $f\underline\#g$, being right linear, corresponds to $x\tl y$,
see \cite{Bm}. In fact, the construction of $A\odot B$ above originates from
properties of the exterior cup products. Compare also the tensor product of
``quadratic modules'' in \cite{CH}.
\end{Rem}

\section{Preliminaries on Nil$_2$-groups}
%\subsection{Nil$_2$-groups}

Groups will be written additively. In particular, for elements $a,b\in G$ of
a group $G$ their commutator will be denoted by $[a,b]=-a-b+a+b$. For any
group $G$ we denote by $\centr(G)$ the center of $G$.
%Moreover we denote by
%$G\ab$ the abelization of $G$, that is, the quotient
%$$G\ab:=$$

A group $G$ is of \emph{nilpotence class two}, or is a
\emph{nil$_2$-group}, if all triple commutators of $G$ vanish, $[[G,G],G]=0$.
The category of all such groups and their homomorphisms will be denoted by
$\Nil$.

For any $G\in\Nil$ there is a well-defined homomorphism $\Lambda^2(G\ab)\to
G$ given by $\hat a\wedge\hat b\mapsto[a,b]$. Here
and elsewhere $\hat x$
denotes the class of $x\in G$ in $G\ab=G/[G,G]$. Moreover, one has the inclusion
$[G,G]\subseteq\centr(G)$ and for any $a,b\in G$ and any $n\in\ZZ$ one has
\begin{equation}\label{nabe}
na+nb=n(a+b)+\binom n2 [a,b].
\end{equation}

The category $\Nil$ has all limits and colimits. For $G_1$ and $G_2$ in
$\Nil$ let $G_1\vee G_2$ denote their coproduct in $\Nil$. Then one has
the following central extension
\begin{equation}\label{jaminilsi}
0\to G_1\ab\t G_2\ab\xto{i} G_1\vee G_2\to G_1\x G_2\to 0.
\end{equation}
Here the homomorphism $i$ is given by $\hat x\t\hat
y\mapsto[i_1(x),i_2(y)]$ for $x\in G_1, y\in G_2$, where
 $i_t:G_t\to G_1\vee G_2$, $t=1,2$
is the canonical inclusion.

The inclusion functor $\Nil\subset\Gr$
has a left adjoint, given by
$$
G\mapsto G\nilp:=G/[[G,G],G].
$$
The forgetful functor $\Nil\to\Sets$ has a left adjoint, whose value on a set
$S$ is known as \emph{the free nilpotent group of class two generated by} $S$
and is denoted by $\brk S\nilp$. One has $\brk S\nilp=(F_S)\nilp$, where
$F_S$ is the free group on $S$.

The following is an easy consequence of the theorem of Witt on the lower
central series of a free group:

\begin{Le}\label{wittori}
For a free nil$_2$-group $G$ one has the central extension
$$
0\to\Lambda^2(G\ab)\to G\to G\ab\to 0.
$$
\end{Le}
%\qed
It follows that there is a normal form of elements in $\brk S\nilp$ for each
linear ordering of the set $S$. Namely, all elements of $\brk S\nilp$ can be
written in a unique way in the form
\begin{equation}\label{nf}
n_1x_1+...+n_px_p+m_1[y_1,z_1]+...+m_q[y_q,z_q]
\end{equation}
with $x_i,y_j,z_j\in S$ and $n_i,m_j\in\ZZ\backslash\set{0}$ for all $i$, $j$
and moreover $x_1<...<x_p$, $y_1<z_1$, ..., $y_q<z_q$ with respect to the
given ordering of $S$ and $(y_1,z_1)<...<(y_q,z_q)$ with respect to the
induced (left) lexicographic ordering of $S\x S$.

%As in any category of algebras,
%monomorphisms in the category $\Nil$ are injective homomorphisms.
%The following lemma shows that epimorphisms in the category $\Nil$ are
%surjective homomorphisms.
%For a morphism $f:G\to H$ in $\Nil$ let $\Cok(f)$ be the quotient of $H$ by the normal subgroup of $H$ generated by $f(H)$.

%\begin{Le}\label{epinil} Let $f:G\to H$ be a homomorphism in $\Nil$. Then the following conditions are equivalent:
%\begin{enumerate}
%\item The homomorphism $f$ is surjective.
%%\item For any morphism $h:H\to H'$ in $\Nil$ with $hf=0$ one has $h=0$.
%
%\item $\Cok(f)=0$.
%\end{enumerate}
%%
%\end{Le}
%\begin{proof} It is clear that $(1)\Longrightarrow (2)\Longrightarrow (3)\Longrightarrow (2)$. Assume the condition $(3)$ holds. Thus for any element $h\in H$ there exists elements
%$h_1,\cdots,h_n\in H$ and $x_1,\cdots,x_n\in G$ such that
%\begin{align*}
%h=&(-h_1+f(x_1)+h_1) +\cdots + (-h_n+f(x_n)+h_n)\\
%=&([h_1,-f(x_1)]+f(x_1))+\cdots +([h_1,-f(x_1)]+f(x_1))\\
%=& [-h_1,f(x_1)]+\cdots +[-h_n,f(x_n)]+f(x_1+\cdots +x_n)
%\end{align*}
%Thus for each $h$ there exists $c\in [H,H]$ and $y\in G$ such that $h=c+f(y)$. In particular we can write $h_i=c_i+f(y_i)$, with $c_i\in [H,H]$ and we obtain
%$$h=[-c_1+f(y_1),f(x_1)]+\cdot +[-c_n+f(y_n),f(x_n)]+f(x_1\cdots +x_n)$$
%Since $c_i$ are central elements, we see that
%$$h=f([y_1,x_1]+\cdots +[y_n,x_n]+x_1+\cdots +x_n).$$
%Thus $f$ is surjective.
%\end{proof}
%A morphism $f:G\to H$ of the category $\Nil$ is called \emph{epimorphism} provided it satisfies the conditions of Lemma \ref{epinil}.

\section{Preliminaries on quadratic maps}\label{4444}
%\subsection{Quadratic maps}\label{4444}

Let $G$ and $G'$ be additively written groups of nilpotence class two.
 Recall that a map $f:G\to G'$ is quadratic if for any
$a,b\in G$ the cross-effect
$$
\cref fab:=-f(b)-f(a)+f(a+b)
$$
is a central element in $G'$ and is linear in $a$ and $b$. Then $f(0)=0$ and the cross-effect yields a
well-defined homomorphism $\cref f--:G\ab\ox G\ab\to \centr(G')$. Moreover, the
following holds (see \cite{niq})
\begin{equation}\label{fminusis}
 f(na)=nf(a)+\binom n2\cref faa , \ \ n\in \ZZ,
 \end{equation}
\begin{equation}\label{fcom}
f([a,b])=\cref fab-\cref fba.
\end{equation}

\begin{Le}\label{extend}
For any set $S$, any abelian group $A$, any map
$$
f_0:S\to A
$$
and any homomorphism
$$
\Ph:\ZZ[S]^{\ox2}\to A
$$
there exists a unique quadratic map
$$
f:\brk S\nilp\to A
$$
satisfying
$$
f(s)=f_0(s)
$$
for $s\in S\subset\brk S\nilp$ and
$$
f(u+v)=f(u)+f(v)+\Ph(\hat u\ox\hat v)
$$
for $u,v\in\brk S\nilp$. In particular there is a unique map
$$
H:\brk S\nilp\to\ZZ[S]\ox\ZZ[S]
$$
satisfying
$$
H(x)=0
$$
for $x\in S$ and
$$
H(u+v)=H(u)+H(v)+\hat v\ox\hat u
$$
for $u,v\in\brk S\nilp$. Thus $H$ is a quadratic map with
\begin{align*}
\cref Huv&=\hat v\ox\hat u,\\
H(-u)&=-H(u)+\hat u\ox\hat u,\\
H([u,v])&=\hat v\ox\hat u-\hat u\ox\hat v.
\end{align*}
\end{Le}

\begin{proof}
Uniqueness is clear from the hypothesis. For the existence, using the above
normal form \eqref{nf}, we can explicitly define
\begin{multline*}
f(n_1x_1+...+n_px_p+m_1[y_1,z_1]+...+m_q[y_q,z_q])\\
=n_1H_0(x_1)+...+n_pH_0(x_p)+\binom{n_1}2\Ph(x_1\ox x_1)+...+\binom{n_p}2\Ph(x_p\ox x_p)+\sum_{1\le i<j\le p}n_in_j
\Ph(x_i\ox x_j)\\
+m_1\Ph(y_1\ox z_1-z_1\ox y_1)+...+m_q\Ph(y_q\ox z_q-z_q\ox y_q).
\end{multline*}
It is easy to see that the so defined $f$ satisfies the required equalities.
\end{proof}

The following key lemma is useful for constructing quadratic maps on groups
given in terms of generators and relations.

\begin{Le}\label{quadex}
Suppose a nil$_2$-group $G$ is given by a set $S=\{x_i;i\in I\}$ of
generators subject to the relations $\{r_j;j\in J\}$. For any abelian group
$A$, any $I$-tuple $(a_i)_{i\in I}$ of elements in $A$ and any homomorphism
$$
\Ph:G\ab\ox G\ab\to A
$$
there exists a unique quadratic map
$$
f:\brk S\nilp
%\ZZ\nilp[(x_i)_{i\in I}]
\to A
$$
satisfying
$$
f(x_i)=a_i,\ \ (x_i|x_{i'})_f=\Ph(x_i,x_{i'}),\ \ i,i'\in I.
$$

Moreover this map factors through the quotient map $q:\brk S\nilp\onto G$
to yield a quadratic map $G\to A$ if and only if it satisfies
$$
f(r_j)=0,\ \ j\in J.
$$
\end{Le}

\begin{proof}
Existence of $f:\brk S\nilp\to A$ is a direct consequence of Lemma
\ref{extend}. It is clear that if $f$ factors trough $G$ then $f(r_j)=0$ for
all $j\in J$. Conversely, assume this condition holds. We have to show that
$f(x)=f(y)$ provided $y=r+x$, where $r$ lies in the smallest normal subgroup
$B$ of $\brk S$ containing all $r_j$, $j\in J$. Observe that $\Phi(r,-)=0$
for any $r\in B$, so that restriction of $f$ to $B$ is a homomorphism and
moreover $f(y)=f(r)+f(x)$, so it remains to prove that $f$ vanishes
identically on $B$. For this, it suffices to show that for any $r\in B$ with
$f(r)=0$ one also has $f(-z+r+z)=0$ for any $z\in\brk S$. Indeed
\begin{align*}
f(-z+r+z)
&=f(-z)+f(r+z)+\Ph(-z,r+z)\\
&=f(-z)+f(z)+\Ph(-z,z)\\
&=f(-z+z)\\
&=0.
\end{align*}
\end{proof}

\section{Square groups and pre-square groups}
In this section we list main properties of square groups and introduce
various special cases of square groups. Also quotients, coproducts, central
extensions in the category of square groups are described. We then outline
the simplicial theory of square groups. Finally we discuss pre-square groups
needed in the proof of our main results.
\subsection{Properties of square groups}\label{squag}
Let
$$
M = (\xymatrix{M\e \ar[r]^H & M\ee \ar[r]^P & M\e})
$$
be a square group as in Section \ref{11}. Then $P\cref Hxy=-x-y+x+y$ implies
that $[M\e,M\e]\subset\Im(P)$, while $\cref H{Pa}y=0=\cref Hx{Pb}$ shows that
$\Im(P)\subset\centr(M\e)$. In particular $M\e$ is a nil$_2$-group and
$\Cok(P)$ is a well-defined abelian group.

\begin{Pro}\label{sqidentities}
\
\begin{itemize}
\item[i)]
One has
$$
H(-x-y+x+y) =\cref Hxy -\cref Hyx.
$$
Moreover, the function
$$
T=HP-\Id
$$
is an involutive automorphism of $M\ee$, i.~e. $T^2=\Id_{M\ee}$. Furthermore
one has
$$
PT=P, \ \ \ T\cref Hxy+\cref Hyx=0.
$$
\item[ii)]
The function $\Delta:M\e\to M\ee$ is linear, where
\begin{align*}
\Delta (x)&=HPH(x)-2H(x)+\cref Hxx\\
&=HPH(x)+H(2x)-4H(x)\\
&=\{x\mid x\}_H-H(x)+TH(x)\\
&=H(-x)+TH(x)
\end{align*}
and furthermore one has
$$
\Delta P = 0,\ \ P\Delta =0, \ \ \Delta+T\Delta=0.
$$
\item[iii)]
For any integer $n$, the map $n^*:M\e \to M\e$ defined by
$n^*(x)=nx+ \binom n2 PH(x)$ is a homomorphism. Moreover one has
$$
(nm)^*=n^*m^*
$$
and
$$
P(n^2a)=n^*(P(a)),\ \ (n^*x\mid n^*y)_H=n^2(x\mid y)_H.
$$
\end{itemize}
\end{Pro}

\begin{proof} i) Since $H$ is quadratic and $M\ee$ is abelian
we can use identity \eqref{fcom} to get the first identity. We have
$$(a\mid b)_T=(a\mid b)_{HP}=(Pa\mid Pb)_H=0.$$
Thus $T$ is a homomorphism. Furthermore, one has
$$T^2=(HP-\Id)(HP-\Id)=HPHP-2HP+\Id=H(P+P)-2HP+\Id=\Id.$$
Similarly
$$PT=P(HP-\Id)=PHP-P=2P-P=P$$
and
$$T\cref Hxy=HP\cref Hxy -\cref Hxy=H([x,y])-\cref Hxy=-\cref Hyx.$$

ii) Since $H$ takes values in an abelian group, we have
\begin{multline*}
\Delta(x+y)=HPH(x+y)-2H(x+y)+\cref H{x+y}{x+y}\\
=HPH(x)+HPH(y)+H([x,y])-2H(x)-2H(y)-2\cref H{x}{y}+\cref H{x+y}{x+y}\\
=\Delta(x)+\Delta(y)+H([x,y])-\cref H{x}{y}+\cref Hyx
=\Delta(x)+\Delta(y).
\end{multline*}
Hence $\Delta$ is additive. To get the other expressions for $\Delta$ observe that $H(2x)=2H(x)+\cref Hxx$ as well as
$TH=HPH-H$ and $H(-x)= -H(x)+\cref Hxx$. Moreover, we have
$\Delta P=HPHP-2HP=0$ and $P\Delta(x)=PHPH(x)-2PH(x)+[x,x]=0$. Similarly
\begin{multline*}T\Delta(x)=THPH(x)-2TH(x)+T(\cref Hxx )=\\
=HPHPH(x)-HPH(x)-2HPH(x)+2H(x)-\cref Hxx\\
=-HPH+2H(x)-\cref Hxx=-\Delta(x).
\end{multline*}
iii) We have

\begin{multline*}
n^*(x+y)=n(x+y)+\binom n2 PH(x+y)\\
=nx+ny+\binom n2(-[x,y]+PH(x)+PH(y)+P(x\mid y)_H)\\
=nx+ny+\binom n2(PH(x)+PH(y))=n^*(x)+n^*(y).
\end{multline*}

Thus $n^*$ is indeed a homomorphism. We also have
$$n^*(m^*(x))= n(mx+\binom m2 PH(x))+\binom n2 PH(mx+\binom m2PH(x)).$$
Since $PH(x)$ is a central element, we obtain

\begin{multline*}
n^*(m^*(x))= nmx+(n\binom m2+\binom n2 m+2\binom n2 \binom m2)PH(x)\\
= nmx +\binom {nm}2 PH(x)=(mn)^*(x).
\end{multline*}

Furthermore, we have
$$n^*(Pa)=nP(a)+\binom n2 PHP(a)=nP(a)+(n^2-n)P(a)=n^2P(a)$$
and

\begin{align*}
(n^*(x)\mid n^*(y))_H&=(nx+P(\binom n2 H(x))\mid ny+P(\binom m2 H(y)))_H\\
&=(nx\mid ny)_H=n^2(x\mid y)_H.
\end{align*}
\end{proof}

\begin{Co}\label{29}
The cross-effect and $\Delta$ yield homomorphisms
$$(-,-)_H:\Cok(P^M)\t\Cok(P^M)\to M\ee$$
and
$$\Delta:\Cok(P^M) \to \Ker(P^M)\subset M\ee.$$
Moreover
$\Delta$ yields the natural homomorphism
$$
k^M:\Cok(P^M)\to \Ker (\frac{M\ee}{\Id-T}\xto{P^M} M\e)
$$
One also has
$$
k^M(\bar x)\equiv(x|x)_H
$$
in $\Ker (\frac{M\ee}{\Id-T}\xto{P^M} M\e)$ for any $x\in M\e$.
\end{Co}
 %%%%%%%%%%%%%%%%%%%%%%%%%%%%%%%%%%%%%%%%%%%%%%%%%%%%%%%%%%%%%%%%%
\subsection{Abelian square groups and quadratic $\ZZ$-modules}\label{qvazimodo}
A square group $M$ is called \emph{abelian} if $H$ is a homomorphism, that is
$\cref Hxy=0$ for all $x,y\in M\e$, equivalently abelian square group
consists of two abelian groups $M\e$ and $M\ee$ together two homomorphisms
$P:M\ee\to M\e$, $H:M\e\to M\ee$ such that $PHP=2P$. Abelian square groups
correspond to quadratic functors $\Gr\to \Ab$ preserving filtered colimits
and reflexive coequalizers. The category of abelian square groups is denoted
by $\Ab(\Square)$.

Let us recall that a quadratic $\ZZ$-module (see \cite{qf}, \cite{square}) is
a square group $M$ for which $\cref H--=0$ and $\Delta=0$. Equivalently a
quadratic $\ZZ$-module is given by two abelian groups $M\e$ and $M\ee$
together with two homomorphisms $P:M\ee\to M_{e}$ $H:M\e\to M\ee$ such that
$PHP=2P$ and $HPH=2H$. Quadratic $\ZZ$-modules correspond to quadratic
functors $\Ab\to \Ab$ preserving filtered colimits and reflexive
coequalizers. The category of quadratic $\ZZ$-modules is denoted by $\QZ$.

Thus we have the following full embeddings:
$$\Ab\subset \QZ\subset \Ab(\Square)\subset \Square.$$
Here abelian groups corresponds to square groups $M$ with $M\ee=0$.

For any abelian group $A$, let $A^\t$ be the quadratic $\ZZ$-module defined
as in Proposition \ref{tlco}.

\begin{Le}\label{eeadj}
For any square group $M$ and for any abelian group $A$ one has the isomorphism
$$\Hom_{\Square}(A^\t,M)\cong \Hom(A,M\ee).$$
\end{Le}
\begin{proof} Take a homomorphism $g:A\to M\ee$. We define
$f=(f\e,f\ee):A^\t\to M$ by $f\e(a)=Pg(a)$ and $f\ee(a,b)=g(a)+Tg(b).$ Then
$f$ is a morphism of square groups and one easily sees that in this way one
gets all such maps. Indeed, one takes $g(a)=f\ee(a,0).$
\end{proof}

Of special interest is the quadratic $\ZZ$-module $\ZZ^\t$ since by Lemma
\ref{eeadj}
$$\Hom_{\sf SQ}(\ZZ^\t,M)\cong M\ee.$$

We will need also the following construction. Let $L$ be an abelian group
and let $\tau$ be an involution on $L$. Then $E(L,\tau)$ is the
quadratic $\ZZ$-module with
$$E(L,\tau)\e=\Cok(L\xto{\Id+\tau} L)$$
$$E(L,\tau)\e=L$$
where $P$ is the natural projection onto quotient, while $H$ is induced by
the homomorphism $\Id-\tau$.

\subsection{Sets versus square groups}\label{lznil}
There is a functor
$$
\Znil[-]:\Sets\to\Square
$$
which is constructed as follows. For a set $S$ one
puts
$$
\Znil[S]\ee=\ZZ[S]\t\ZZ[S],
$$
where $\ZZ[S]$ is the free abelian group generated by $S$. We
take $\Znil[S]\e$ to be $\brk S\nilp$, the free nil$_2$-group generated by $S$.
The homomorphism $P$ is given by $P(s\t t)=[t,s]$, $s,t\in S$, while the quadratic map
$H$ is uniquely defined by
$$
H(s)=0, \ \ (s\mid t)_H=t\t s\ \ s,t\in S.
$$
If $S$ is a singleton, we obtain $\Znil=\Znil[S]$, see Definition \ref{sgdef}.

For general $S$ one has
$$
\Cok(P^{\Znil[S]})=\ZZ[S]
$$
and the homomorphism
$$
\Delta:\ZZ[S]=\Cok(P^{\Znil[S]})\to\ZZ[S]\t\ZZ[S]=\Znil[S]\ee
$$
is given by $\Delta(s)=(s,s)$. Moreover, the homomorphism
$$
T:\ZZ[S]\t\ZZ[S]\to\ZZ[S]\t\ZZ[S]
$$
is given by $T(s\t t)=-t\t s$.
%%%%%%%%%%%%%%%%%%%%%%%%%%%%%%%%%%%%%%%%%%%%%%%%%%%%%%%%%%%%
%%%%%%%%%%%%%%%%%%%%%%%%%%

It turns out that the functor $\Znil[-]:\Sets\to \Square$ is a left adjoint.
 Let $M$ be a square group. An element
$x\in M\e$ is called \emph{linear} if $H(x)=0$. Let $\LL(M)$ be the subset of
linear elements in $M\e$ so that we obtain a functor
$$
\LL:\Square\to\Sets.
$$

\begin{Pro}\label{slineareleents}
The functor $\Znil[-]:\Sets\to \Square$ is left adjoint to the functor
$\LL$.
\end{Pro}

\begin{proof}
Let $S$ be a set and $M$ be a square group. Given a morphism $f:\Znil[S]\to
M$ of square groups, the composite $S\subset\Znil[S]\to M$ is a map $f_0$;
moreover its image is contained in $\LL(M)$ since $H(s)=0$ for $s\in S$ and
$f$ is compatible with $H$.

Conversely we must show that any map $f_0:S\to\LL(M)\subset M\e$ extends
uniquely to a square group morphism $f:\Znil[S]\to M$. First, there is
clearly a unique group homomorphism $f\e:\Znil[S]\e\to M\e$ extending $f_0$,
as $\Znil[S]\e$ is the free nil$_2$-group on $S$ and $M\e$ is a
nil$_2$-group. Moreover, by compatibility of a morphism of square groups with
$H$ we necessarily have
$$
f\ee(s\ox s')=f\ee((s'|s)_H)=(f_0(s')|f_0(s))_H
$$
for any $s,s'\in S$. Hence we also have a unique choice for
$f\ee:\ZZ[M]\ox\ZZ[M]\to M\ee$ and one has
$$
f\ee(\bar x\ox\bar y)=(f\e(y)|f\e(x))_H
$$
for any $x,y\in\Znil[M]\e$.

It remains to show that the pair $(f\e,f\ee)$ yields a
morphism of square groups
$$
(f\e,f\ee):\Znil[S]\to M.
$$
Indeed, compatibility with $P$ is clear since
$$
f\e P(s\ox s')=f\e P((s'|s)_H)=f\e[s',s]=[f_0(s'),f_0(s)]=P((f_0(s')|f_0(s))_H)%\\
=Pf\ee(s\ox s')
$$
as $f\e$, $f\ee$ and $P$ are group homomorphisms. Since image of $f_0$ is in
$\LL(M)$, compatibility with $H$ holds on elements of $S$; moreover if it
holds on $x$ and $y$, one has
\begin{multline*}
f\ee H(x+y)=f\ee(H(x)+H(y)+\bar y\ox\bar x)=Hf\e(x)+Hf\e(y)+(f\e(x)|f\e(y))_H\\
=H(f\e(x)+f\e(y))=Hf\e(x+y).
\end{multline*}
This finishes the proof.
\end{proof}

\subsection{Normal subobjects and quotients of square groups}\label{quot}
Obviously for any morphism $f:M\to N$ of square groups kernels of $f\e$ and
$f\ee$ determine a sub-square group $\Ker(f)$ of $M$. Sub-square groups of
this form can be characterized as those $K\into M$ for which $K\e$ is normal
in $M\e$ and moreover one has
$$
(M\e|K\e)_H, (K\e\mid M\e)_H \subset K\ee,
$$
i.~e. for any $x\in M\e$, $k\in K\e$ one has $(x|k)_H\in K\ee$
and $(k|x)_H\in K\ee$. Such
sub-square groups will be called \emph{normal}. For any normal sub-square
group $K\normal M$ the quotient $M/K$ is defined, with $(M/K)\e=M\e/K\e$,
$(M/K)\ee=M\ee/K\ee$; here $P:M\ee/K\ee\to M\e/K\e$ is the uniquely determined
homomorphism whereas $H:M\e/K\e\to M\ee/K\ee$ is the uniquely determined map
by virtue of
$$
H(x+k)-H(x)=H(k)+(x|k)_H\in K\ee
$$
for any $x\in M\e$, $k\in K\e$.

Cokernels of morphisms in $\Square$ are defined as follows. For
a morphism $f:M\to N$ of square groups, let $\Cok(f)$ be the quotient of $N$ by
the smallest normal sub-square group generated by $\Im(f)$. Thus one has
$$
\Cok(f)\e=\Cok(f\e),
$$
i.~e. $\Cok(f)\e$ is the quotient of $N\e$ by the normal subgroup generated
by the image of $f\e$, and
$$
\Cok(f)\ee=\Cok(f\ee)/(\Im(f\e)|N)_H,
$$
that is, $\Cok(f)\ee$ is the quotient of $N\ee$ by elements of the form
$f\ee(a)$ for $a\in M\ee$ and $(f\e(x)|y)_H$ for $x\in M\e$, $y\in N\e$.

\begin{Le}\label{episq}
Let $f:M\to N$ be a morphism in $\Square$. Then the following conditions are
equivalent:
\begin{enumerate}
\item
$f\e$ and $f\ee$ are surjective;
\item
for any morphism $h:N\to N'$ in $\Square$ with $hf=0$ one has $h=0$;
\item
$\Cok(f)=0$.
\end{enumerate}
\end{Le}

\begin{proof}
We show that $(3)\Longrightarrow (1)$. The rest is trivial.
Assume (3) holds. It follows that $\Cok(f\e)=0$. Hence $f\e$
is surjective (see for example Exercise 5, Section 5, Chapter 1 in \cite{working}).
Moreover (3) implies $$N\ee=f\ee(M\ee) + (f\e(M\e)\mid N\e)_H
+(N\e\mid f\e(M\e))_H.$$ Since $N\e=f\e(M\e)$, we see that
$$
N\ee=f\ee(M\ee)+\cref H{f\e(M\e)}{f\e(M\e)}=f\ee(M\ee)+f\ee(\cref
H{M\e}{M\e})
$$
and hence $f\ee$ is surjective.
\end{proof}

A morphism $f:M\to N$ in the category $\Square$ is called an
\emph{epimorphism} provided it satisfies the conditions of Lemma \ref{episq}.
One easily deduces from Lemma \ref{episq} that the class of effective
epimorphisms (Section 4, Chapter 2, \cite{qu}) in the category $\Square$
coincides with the class of epimorphisms.

\subsection{Central extensions of square groups}
A sequence of square groups
$$0\to A\to B\to C\to 0$$
is called \emph{short exact} if it is exact on the $\mathrm{e}$-level and the
$\mathrm{ee}$-level. We will say that it is a \emph{central extension} if
$A\e$ is a central in $B\e$ and $\cref Hxy=0$ provided $x\in A\e$ and $y\in
B\e$, or $x\in B\e$ and $y\in A\e$. In particular $A$ is a normal sub-square
group of $B$ and $H$ is linear on $A\e$.
\begin{Le}\label{centralizator} Let
$$0\to A\to B\to C\to 0$$
be a short exact sequence in $\SG$. Then there is a well-defined
 square group $A'$ defined by
$$A'\e:=\{x\in A\e\mid x+y=y+x\ \& \ \cref Hxy=0=\cref Hyx, \ y\in B\e\},$$
$$A'\ee:=\{a\in A\ee\mid P(a)\in A'\e\},$$
and $H$ and $P$ being the restriction of $P^A$ and $H^A$. Moreover the
columns and the bottom row of the commutative diagram
$$\xymatrix{&0\ar[d]&0\ar[d]&&\\
&A'\ar[r]^{\Id}\ar[d]&A'\ar[d]&&\\
0\ar[r]&A\ar[r]\ar[d]&B\ar[r]\ar[d]&C\ar[r]\ar[d]^{\Id}&0\\
0\ar[r]&A/A'\ar[r]\ar[d]&B/A'\ar[r]\ar[d]&C\ar[r]&0\\
&0&0&&
}$$
are central extensions of square groups.
\end{Le}

\begin{proof}
Take any element $x$ from $A\e$. Then $PH(x)\in A\e'$. In particular
$H(A\e')\subset A\ee'$ and therefore $A'$ is well-defined. It is obvious that
the columns are central extensions of square groups. Take now any elements
$x\in A\e$ and $y\in B\e$. Then the commutator $[x,y]$ projects to zero in
$C$ and therefore $[x,y]\in A\e$. Since $N\e\in\Nil$ and the cross-effect of
$H$ vanishes on commutators it follows that $[x,y]\in A\e'$. Thus $(A/A')\e$
is a central subgroup of $(B/A')\e$. It remains to show that $H(x,y)\in
A'\ee$. But this follows immediately from the facts that the image of $H(x,y)$
in $C\ee$ vanishes (thus $H(x,y)\in A\ee$) and $PH(x,y)\in A\e'$.
\end{proof}

\subsection{Coproduct of square groups}\label{copr}
 Let $M$ and $N$
be square groups. Then their coproduct $M\vee N$ in
the category of square groups has the
following form
$$M\vee N=((M\vee N)\e\xto{H} (M\vee N)\ee\xto{P}
(M\vee N)\e)$$
where
$$(M\vee N)\ee=M\ee\oplus N\ee\oplus \Cok(P^M)\t \Cok(P^N)\oplus \Cok(P^N)\t
\Cok(P^M)$$
while $(M\vee N)\e$ is the quotient of the coproduct $M\e\vee N\e$ in the
category $\Nil$ by the following relations

\begin{align*}
P^M(a)+y&=y+P^M(a), &a\in M\ee, y\in N\e,\\
P^N(b)+x&=x+P^N(b), &b\in N\ee, x\in M\e.
\end{align*}
Moreover, $P$ and $H$ of $M\vee N$ are given by
\begin{align*}
P(a+b+x_1\t y_1+y_2\t x_2)&=P^M(a)+P^N(b)+[x_1,y_1]+[y_2,x_2],\\
H(x+y+[x_1,y_1])&=H^M(x)+H^N(y)+x\t y+x_1\t y_1-y_1\t x_1
\end{align*}
where $x,x_1,x_2\in M\e$, $y,y_1,y_2\in N\e$, $a\in M\ee$, $b\in N\ee$.

Let us also observe that one has a central extension of the form
$$
0\to \Cok(P^M)\t \Cok(P^N)\to (M\vee N)\e\to M\e\x N\e\to 0,
$$
which implies the short exact sequence of square groups
\begin{equation}\label{sqcopr}
0\to (\Cok(P^M)\t \Cok(P^N))^\t \xto{j} M\vee N\to M\x N\to 0.
\end{equation}
Here $j\e(x\t y)=[x,y]$ and $j\ee(x_1\t y_1,x_2\t y_2)=x_1\t y_1+y_2\t x_2$.

Since the map $P^M$ is surjective for $M=A^\t$, we obtain that for any
abelian group $A$ and any square group $M$ one has
$$A^\t\vee M \cong A^\t \x M.$$
In particular for abelian groups $A$ and $B$ one gets
$$A^\t\vee B^\t \cong (A\oplus B)^\t.$$

\subsection{Free and projective square groups}\label{Vconst}

We need the following square group $\ZZ^Q$ defined by
 $$(\ZZ^Q)\e= \ZZ\oplus \ZZ, \ \
(\ZZ^Q)\ee=\ZZ\oplus \ZZ\oplus \ZZ,$$ with the maps $P$ and $H$ given by
$P(a,b,c)=(0,a+2b)$ and $H(m,n)=(m,n,\binom m2)$. One easily shows that for a
square group $M$ one has the natural isomorphism
$$\Hom_{\Square}(\ZZ^Q,M)\cong M\e.$$
For any set $S$ we put
$$V(S):=\bigvee_{s\in S}\ZZ^Q.$$
It follows that
$$
\Hom_{\Square}(V(S),M) \cong \Hom_{\sf Sets}(S,M\e).
$$
Thus the functor $V:{\sf Sets}\to \Square$ is left adjoint to the functor
$$
\Square \to {\sf Sets}, \ \ M\mapsto M\e.
$$
Now we give the following explicate construction of $V(S)$. We consider three
further copies of $S$, which are denoted respectively by $HS$, $PHS$ and
$HPHS$. For an element $s\in S$ the elements $Hs$, $PHs$, $HPHs$ correspond
to $s$ in these copies. Then we take
\begin{align*}
V(S)\e&=\brk S\nilp\x \ZZ[PHS],\\
V(S)\ee&=\ZZ[HS]\oplus\ZZ[S\x S]\oplus \ZZ[HPHS].
\end{align*}
Moreover, $H$ is the unique quadratic map with
\begin{align*}
H(s)&=Hs,\\
H(PHs)&=HPHs,\\
\cref Hst&=(s,t),\\
\cref Hs{PHt}&=0=\cref H{PHs}t=\cref H{PHs}{PHt}.
\end{align*}
Here $s,t\in S$. The homomorphism $P$ is given by
$$
P(Hs)=PHs, \ \ P(s,t)=[s,t], \ \ P(HPHs)=2PHs.
$$
The fact that this is really isomorphic to $V(S)$ can be deduced either from
Section \ref{copr} or directly by the universal property.

Properties of the functor $A\mapsto A^\t$ imply that the functor ${\sf Sets\x
Sets}\to \Square$ given by
$$
(S,T)\mapsto V(S)\x (\ZZ[T])^\t\cong V(S)\vee (\ZZ[T])^\t
$$
is left adjoint to the forgetful functor
$$
\Square \to {\sf Sets\x Sets}, \ \ M\mapsto (M\e,M\ee).
$$
A square group is called \emph{free} if it is isomorphic to $V(S)\x
(\ZZ[T])^\t$. A square group is called \emph{projective} provided it has the
familiar lifting property with respect to epimorphisms of square groups. Any
free square group is projective and any projective square group is a retract
of a free square group.
\begin{Le}
For any square group $M$ there exists an epimorphism $F\onto M$, where $F$ is
free.
\end{Le}

\begin{proof} One can take $F=V(M\e)\x (\ZZ[M\ee])^\t$ with the morphism
adjoint to $(\Id_{M\e},\Id_{M\ee})$.
\end{proof}

\subsection{Simplicial objects in the category $\Square$}\label{ssg}
Let $\ssq$ be the category of simplicial objects in the category $\Square$ of
square groups. Any such simplicial object $X$ defines two simplicial groups
$X\e$ and $X\ee$ as well as a morphism $P:X\ee\to X\e$ of simplicial groups.
The map $H$ yields a morphism of simplicial sets $H: X\e\to X\ee$. If one
passes to homotopy groups, then one obtains groups $\pi_i(X\e)$ and
$\pi_i(X\ee)$ together with induced homomorphisms $P:\pi_i(X\ee)\to
\pi_i(X\e)$, $i\ge0$. The map $H$ yields a quadratic map $H:\pi_0(X\e)\to
\pi_0(X\ee)$ and homomorphisms $H:\pi_i(X\e)\to\pi_i(X\ee)$ for $i\ge1$. It
is clear that the equation $PHP=2P$ still holds for induced maps. It follows
that for each $i\ge1$ one obtains a well-defined abelian square group
$$
\pi_iX\in\Ab(\Square), i\ge1
$$
with
$$
(\pi_iX )\e=\pi_i(X\e), \ \ (\pi_iX )\ee=\pi_i(X\ee).
$$
We also have a well-defined square group
$$
\pi_0(X)\in \Square,
$$
since all equations defining a square group hold in $X_0$, the zero component
of $X$, and therefore they remain true in the quotient $\pi_0(X)$.

By \ref{Vconst} the category $\Square$ satisfies all conditions of Theorem 4
Section 4, Chapter 2 \cite{qu}. Hence the category $\ssq$ of simplicial
objects of $\Square$ possesses a closed model category structure, where a
morphism $f$ is a weak equivalence (resp. fibration) provided $f\e$ and
$f\ee$ are weak equivalences (resp. fibrations) of underlying simplicial
sets. According to \cite{qu} cofibrations are retracts of free maps as they
are defined in \cite{qu}. Equivalently, a map $f:X\to Y$ is weak equivalence
if and only if the induced map of square groups $f_i:\pi_i(X)\to \pi_i(Y)$ is
an isomorphism for all $i\ge0$.

\subsection{Pre-square groups}
\begin{De}\label{psg} A \emph{pre-square group} consists of a diagram
$$
M=\left(M\e\times M\e\xto{\{-,-\}}M\ee\xto TM\ee\xto PM\e\right).
$$
Here $M\ee$ is an abelian group and $T$ is a homomorphism with $T^2=\Id$ and
 $M\e$ is a group written additively, $P$ is a homomorphism and
$\{-,-\}$ is a bilinear map, that is $\{x+y,z\}=\{x,z\}+\{y,z\}$ and
$\{x,y+z\}=\{x,y\}+\{x,z\}$, for all $x,y,z\in M\e$. Moreover one requires
the following identities:
\begin{enumerate}
\item[(a)] $PT=P$,
\item[(b)] $T\{x,y\}+\{y,x\}=0$, $x,y\in M\e$,
\item[(c)] $P\{x,y\}=-x-y+x+y$, $x,y\in M\e$,
\item[(d)] $\{x,Pa\}=0$, $x\in M\e$, $a\in M\ee$.
\end{enumerate}
\end{De}
It follows from (b) that one has $\{Pa,x\}=0$. It follows from (c) and (d)
that $Pa$ lies in the center of $M\e$. Thus $\Cok(P)$ is well-defined and by
(c) it is an abelian group. It follows that $M\e$ is a group of nilpotence
class 2. Bilinearity of the bracket together with (d) shows that there is a
well-defined homomorphism
$$
\{-,-\}:\Cok(P)\t \Cok(P)\to M\ee.
$$

We let $\PSG$ denote the category of pre-square groups. It is clear that the
full subcategory of $\PSG$ consisting of pre-square groups with trivial
$M\ee=0$ is equivalent to the category of abelian groups. In what follows we
identify abelian groups with such pre-square groups.

Thanks to Proposition \ref{sqidentities}, for any square group the following object
$$
\wp(M)=(M\e,M\ee,T=HP-\Id,(-,-)_H,P)
$$
is a pre-square group. Thus we obtain the forgetful functor
$$
\wp:\Square \to \PSG.
$$

Comparing the definitions we immediately obtain the following easy, but useful
result.

\begin{Le}\label{pre=square}
Let $M$ be a pre-square group and let $H:M\e\to M\ee$ be a map.
Then $(M\e,M\ee, P, H)$ is a square group with $\wp(M\e,M\ee,P,H)=M$
iff $\{x,y\}=\cref Hxy$ and $\Id+T=HP$.
\end{Le}

The following Lemma is an immediate consequence of Proposition \ref{sqidentities}.

\begin{Le}\label{intpsg}
For any square group $M$ and any integer $n$, there is a morphism
of pre-square groups
$$
n^*=(n^*)_M:\wp(M)\to\wp(M)
$$
which on $\mathrm{ee}$-level is the multiplication by $n^2$, while on
$\mathrm e$-level it is given by
$$
x\mapsto nx+\binom n2 PH(x).
$$
\end{Le}

\section{The tensor product as a monoidal structure of $\Square$}\label{tl901}
In this section we prove Theorem \ref{tlmonoidal}.

\subsection{The tensor product of a pre-square group and a square group}\label{tlpsg}
Before considering the tensor product of square groups as a square group for
technical reasons
we construct first the bifunctor $$\tl:\PSG\x \Square \to \PSG.$$
\begin{De}\label{tlactiond}
Let $M$ be a pre-square group and $N$ be a square group, then $M\tl N$ is the
pre-square group, which on the $\mathrm{ee}$-level is given by
$$
(M\tl N)\ee=M\ee\t N\ee.
$$
Moreover $(M\tl N)\e$ is generated by elements of the form $x\tl y$ for $x\in
M\e$, $y\in N\e$ and $a\tb b$ for $a\in M\ee$, $b\in N\ee$, subject to the
relations
\begin{enumerate}
\item\label{bice} the symbol $a\tb b$ is bilinear and central in $(M\tl N)\e$,
\item\label{ledi} $x\tl(y_1+y_2)=x\tl y_1+x\tl y_2$,
\item\label{ridi} $(x_1+x_2)\tl y=x_1\tl y+x_2\tl y+\{x_2, x_1\}\tb H(y)$,
\item\label{rip} $x\tl P(b)=\{x,x\}\tb b$,
\item\label{lep} $P(a)\tl y=a\tb\Delta(y)$,
\item\label{TT} $T(a)\tb T(b)=-a\tb b$.
\end{enumerate}
The homomorphism
$$
P:(M\tl N)\ee\to(M\tl N)\e
$$
is given by
$$
P(a\ox b)=a\tb b.
$$
The involution on $(M\tl N)\ee$ is given by
$$
T(a\t b)=-T(a)\t T(b),
$$
while the bracket is given by
$$
\{x\tl y,x'\tl y'\}=\{x,x'\}\t (y\mid y')_H,
$$
$$
\{u,a\tb b\}=0=\{a\tb b,u\},
$$
where $u\in (M\tl N)\e$.
\end{De}

We now show that $M\tl N$ is a well-defined pre-square group satisfying the
identities in Definition \ref{psg}. Identity a) follows from the identity 6)
of the definition of $M\tl N$, while d) is a direct consequence of the
description of $P$. Identity b) can be verified as follows
\begin{align*}
T\{x\tl y,x'\tl y'\}=& T(\{x,x'\}\t (y\mid y')_H)\\
=&-T(\{x,x'\})\t T((y\mid y')_H)\\
=&-\{x',x\}\t (y'\mid y)_H\\
=&-\{x'\t y',x\t y\}.
\end{align*}
It remains to check identity c). Since $P\{x\tl y,x'\tl y'\}=\{x,x'\}\tb
(y\mid y')_H$, this must be equal to $-x\tl y-x'\tl y'+x\tl y+x'\tl y'$.
Consider
\begin{align*}
(x'+x)\tl(y+y')=&(x'+x)\tl y+(x'+x)\tl y'\\
=&x'\tl y+x\tl y+\{x, x'\}\tb H(y)+x'\tl y'+x\tl y'+\{x,x'\}\tb H(y')\\
=&x'\tl y+x\tl y+x'\tl y'+x\tl y'+\{x,x'\}\tb(H(y)+H(y'))\\
=&x'\tl y+x\tl y+x'\tl y'+x\tl y'
+\{x,x'\}\tb(H(y+y'))-\{x,x'\}\tb(y|y')_H.
\end{align*}
On the other hand the same expression expands to
\begin{align*}
(x'+x)\tl(y+y')&=x'\tl(y+y')+x\tl(y+y')+\{x,x'\}\tb H(y+y')\\
&=x'\tl y+x'\tl y'+x\tl y+x\tl y'+\{x,x'\}\tb H(y+y').
\end{align*}
Comparing these expressions gives
\begin{equation}\label{tlcom}
-x\tl y-x'\tl y'+x\tl y+x'\tl y'=\{x,x'\}\tb(y|y')_H
\end{equation}
which is the equality we need. Thus we have constructed a well-defined tensor
product $\PSG\x\Square\to\PSG.$

\begin{Le}\label{tenzigiveoba} Let $M$ be a pre-square group and $N$
be a square group. Then one has the following identities in $M\tl N$.
\begin{enumerate}
\item $\{x_2,x_1\} \tb b=\{x_1,x_2\} \tb Tb$,
\item $(nx)\tl y=x\tl (ny+\binom n2 PHy)$,
\item $[x_1,x_2]\tl y=\{x_1,x_2\}\tb (HPHy-2Hy).$
\end{enumerate}
Here $x,x_1,x_2\in M\e$, $y\in N\e$, $b\in M\ee$ and $n\in \ZZ$.
\end{Le}

\begin{proof} We have
$$\{x_2,x_1\} \tb b=(-T\{x_1,x_2\})\tb TTb=\{x_1,x_2\} \tb Tb$$
and (1) is proved. For a given $y\in N\e$
consider the map $f:M\e\to M\ee$ defined by $f(x)=x\tl y$. Then $f$ is
quadratic, with cross-effect given by
$$\cref f{x_1}{x_2} =\{x_2,x_1\}\tb H(y).$$
Thus by identity \eqref{fminusis} in Section \ref{4444} we have
$$(nx)\tl y=f(nx)=nf(x)+\binom n2 \cref fxx=n(x\tl y)+
\binom n2(\{x,x\})\tb H(y).$$ Since $\tl$ is linear with respect to the
second variable, (2) follows. Next we use identity \eqref{fcom} in section
\ref{4444} to get
$$[x_1,x_2]\tl y=f([x_1,x_2])=\cref f{x_1}{x_2}-\cref f{x_2}{x_1}=\{x_2,x_1\}\tb Hy-
\{x_1,x_2\}\tb Hy.$$
By (1) in Lemma \ref{tenzigiveoba} we get
$$[x_1,x_2]\tl y=\{x_1,x_2\}\tb (THy-Hy).$$
But $TH=HPH-H$ and the result follows.
\end{proof}

\subsection{The tensor product of square groups}\label{hisarseboba}

Assume now that $M$ and $N$ are square groups. We have to show that $H_{M\tl
N}$ is well-defined and $M\tl N$ is in fact a square group. To prove the
first assertion we show that the conditions of Lemma \ref{quadex} are indeed
satisfied, namely, that the quadratic map $H$ (when considered as a quadratic
map from a free nil$_2$-group to $(M\ox N)\ee$) respects the relations from
the definition. First, bilinearity of $a\tb b$ is respected since $H(a\tb
b)=a\ox b-T(a)\ox T(b)$ is bilinear and the cross-effect of $H$ vanishes on
all elements of the form $a\tb b$. Centrality of $a\tb b$ is trivially
respected as the values are taken in an abelian group. Next, the relation
\eqref{ledi} in Definition \ref{tldef} is respected since
\begin{align*}
H(x\tl y_1&+x\tl y_2-x\tl(y_1+y_2))\\
=&H(x\tl y_1+x\tl y_2)+H(-x\tl(y_1+y_2))+(x\tl y_1+x\tl y_2|-x\tl(y_1+y_2))_H\\
=&H(x\tl y_1)+H(x\tl y_2)+(x\tl y_1|x\tl y_2)_H-H(x\tl(y_1+y_2))\\
&+(x\tl(y_1+y_2)|x\tl(y_1+y_2))_H+(x\tl y_1+x\tl y_2|-x\tl(y_1+y_2))_H\\
=&(x|x)_H\ox H(y_1)+H(x)\ox\Delta(y_1)+(x|x)_H\ox H(y_2)+H(x)\ox\Delta(y_2)\\
&+(x|x)_H\ox(y_1|y_2)_H-(x|x)_H\ox H(y_1+y_2)\\
&-H(x)\ox\Delta(y_1+y_2)+(x|x)_H\ox(y_1+y_2|y_1+y_2)_H\\
&-(x\tl y_1|x\tl(y_1+y_2))_H-(x\tl y_2|x\tl(y_1+y_2))_H\\
=&(x|x)_H\ox H(y_1)+(x|x)_H\ox H(y_2)+(x|x)_H\ox(y_1|y_2)_H\\
&-(x|x)_H\ox(H(y_1)+H(y_2)+(y_1|y_2)_H)\\
&+(x|x)_H\ox(y_1|y_1)_H+(x|x)_H\ox(y_1|y_2)_H\\
&+(x|x)_H\ox(y_2|y_1)_H+(x|x)_H\ox(y_2|y_2)_H\\
&-(x|x)_H\ox(y_1|y_1+y_2)_H-(x|x)_H\ox(y_2|y_1+y_2)_H\\
&=0.
\end{align*}

For the relation \eqref{ridi} Definition \ref{tldef} we have
\begin{align*}
H(x_1\tl y&+x_2\tl y+(x_2|x_1)_H\tb H(y)-(x_1+x_2)\tl y)\\
=&H(x_1\tl y+x_2\tl y+(x_2|x_1)_H\tb H(y))+H(-(x_1+x_2)\tl y)\\
&+(x_1\tl y+x_2\tl y+(x_2|x_1)_H\tb H(y)|-(x_1+x_2)\tl y)_H\\
=&H(x_1\tl y+x_2\tl y)+H((x_2|x_1)_H\tb H(y))\\
&-H((x_1+x_2)\tl y)+((x_1+x_2)\tl y|(x_1+x_2)\tl y)_H\\
&-(x_1\tl y|(x_1+x_2)\tl y)_H-(x_2\tl y|(x_1+x_2)\tl y)_H\\
=&H(x_1\tl y)+H(x_2\tl y)+(x_1\tl y|x_2\tl y)_H\\
&+(x_2|x_1)_H\ox H(y)-T((x_2|x_1)_H)\ox TH(y)\\
&-(x_1+x_2|x_1+x_2)_H\ox H(y)-H(x_1+x_2)\ox\Delta(y)\\
&+(x_1+x_2|x_1+x_2)_H\ox(y|y)_H\\
&-(x_1|x_1+x_2)_H\ox(y|y)_H-(x_2|x_1+x_2)_H\ox(y|y)_H\\
=&(x_1|x_1)_H\ox H(y)+H(x_1)\ox\Delta(y)+(x_2|x_2)_H\ox H(y)+H(x_2)\ox\Delta(y)\\
&+(x_1|x_2)_H\ox(y|y)_H+(x_2|x_1)_H\ox H(y)+(x_1|x_2)_H\ox TH(y)\\
&-(x_1+x_2|x_1+x_2)_H\ox H(y)-(H(x_1)+H(x_2)+(x_1|x_2)_H)\ox\Delta(y)\\
=&(x_1|x_2)_H\ox(y|y)_H+(x_1|x_2)_H\ox TH(y)\\
&-(x_1|x_2)_H\ox H(y)-(x_1|x_2)_H\ox\Delta(y)\\
=&0.
\end{align*}

Next for the relation \eqref{rip} Definition \ref{tldef} we check
\begin{multline*}
H(x\tl P(b)-(x|x)_H\tb b)=H(x\tl P(b))-H((x|x)_H\tb b)\\
=(x|x)_H\ox HP(b)+H(x)\ox\Delta P(b)-(x|x)_H\ox b+T((x|x)_H)\ox T(b)\\
=(x|x)_H\ox(HP(b)-b-T(b))=0
\end{multline*}
and for \eqref{lep}
\begin{multline*}
H(P(a)\tl y-a\tb\Delta(y))=H(P(a)\tl y)-H(a\tb\Delta(y))\\
=(P(a)|P(a))_H\ox H(y)+HP(a)\ox\Delta(y)-a\ox\Delta(y)+T(a)\ox T\Delta(y)\\
=(HP(a)-a-T(a))\ox\Delta(y)=0.
\end{multline*}

Finally the relation \eqref{TT} Definition \ref{tldef} is respected since
\begin{multline*}
H(a\tb b+T(a)\tb T(b))=H(a\tb b)+H(T(a)\tb T(b))\\
=a\ox b-T(a)\ox T(b)+T(a)\ox T(b)-TT(a)\ox TT(b)=0.
\end{multline*}

Moreover we have to show that identities of square groups hold for $M\tl N$.
But we have already proved that it is a pre-square group, thus by Lemma
\ref{pre=square} we have only to check the identity $\Id+T=HP$, which holds
because
$$
HP(a\ox b)=H(a\tb b)=a\ox b-T(a)\ox T(b))=a\tb b+T(a\ox b).
$$
%\end{proof}

Thus we get the well-defined tensor product $-\tl-:\Square \x \Square\to\Square.$

\begin{Le}\label{kidevtenzigiveoba} Let $M$ and $N$ be
square groups. Then one has the following identity in $M\tl N$.
$$HPH(x)\tb b=H(x)\tb (b-Tb).$$
Here $x\in M\e$ and $b\in M\ee$.
\end{Le}

\begin{proof} We have
$$(HPHx)\tb b=(Hx+THx)\tb b=(Hx)\tb b+(THx)\tb b=(Hx)\tb b-(Hx)\tb Tb$$
and the result follows.

\end{proof}

\subsection{Associativity}

Our next goal is to show that the tensor product in $\Square$ defines a
symmetric monoidal structure. To construct associativity
isomorphisms we introduce the triple
tensor product $A\tl B\tl C$; we will then construct isomorphisms of this
object to $(A\tl B)\tl C$ and $A\tl(B\tl C)$.

We define $(A\tl B\tl C)\e$ by generators of the form $x\tl y\tl z$ for $x\in
A\e$, $y\in B\e$, $z\in C\e$ and $a\tb b\tb c$ for $a\in A$, $b\in B$, $c\in
C$, subject to the relations
\begin{enumerate}
\item $a\tb b\tb c$ is central and trilinear;
\item $x\tl y\tl(z+z')=x\tl y\tl z+x\tl y\tl z'$;
\item $x\tl(y+y')\tl z=x\tl y\tl z+x\tl y'\tl z+(x|x)_H\tb(y'|y)_H\tb H(z)$;
\item $(x+x')\tl y\tl z=x\tl y\tl z+x'\tl y\tl z+(x'|x)_H\tb(y|y)_H\tb H(z)+(x'|x)_H\tb H(y)\tb\Delta(z)$;
\item $P(a)\tl y\tl z=a\tb\Delta(y)\tb\Delta(z)$;
\item $x\tl P(b)\tl z=(x|x)_H\tb b\tb\Delta(z)$;
\item $x\tl y\tl P(c)=(x|x)_H\tb(y|y)_H\tb c$;
\item $T(a)\tb T(b)\tb T(c)=a\tb b\tb c$.
\end{enumerate}
Moreover we define $(A\tl B\tl C)\ee=A\ox B\ox C$ and $P(a\ox b\ox c)=a\tb
b\tb c$. Finally we let $H$ be the unique quadratic map satisfying
$$
H(x\tl y\tl z)=
(x|x)_H\tb(y|y)_H\tb H(z)+(x|x)_H\tb H(y)\tb\Delta(z)+H(x)\tb\Delta(y)\Delta(z)
$$
with cross-effect equal to
$$
\rho: \Cok(P_{A\tl B\tl C})\t \Cok(P_{A\tl B\tl C})
\xto\rho A\ee\ox B\ee\ox C\ee
$$
given by
$\rho(\bar a_1\ox\bar b_1\ox\bar c_1 \t \bar a_2\ox\bar b_2\ox\bar c_2)=
(a_1|a_2)_H\t (b_1|b_2)_H\t (c_1|c_2)_H$. Here we use the identification
\begin{multline*}
\Cok(P_{A\tl B\tl C})\t \Cok(P_{A\tl B\tl C})\\
=\Cok(P_A)\t \Cok(P_B)\t \Cok(P_C)\t \Cok(P_A)\t
\Cok(P_B)\t \Cok(P_C)
\end{multline*}

Again by Lemma \ref{quadex} such $H$ exists and is unique and the same
argument as in Section \ref{hisarseboba} shows that $A\tl B\tl C$ is a
well-defined square group. By the same methods one shows that there is
the unique morphism of square groups
$$
\alpha=\alpha_{A,B,C}: ((A\tl B)\tl C)\to A\tl B\tl C
$$
which is the canonical isomorphism $(A\ee\t B\ee)\t C\ee\to A\ee\t B\ee\t
C\ee$ on the ee-level, while on the e-level it satisfies the identities
\begin{align*}
\alpha((x\tl y)\tl z)&=x\tl y\tl z,\\
\alpha((a\tb b)\tl z)&=a\tb b \tb \Delta(z),\\
\alpha( (a\t b)\tb c)&=a\tb b\tb c.
\end{align*}
This is an isomorphism with inverse given by
\begin{align*}
\alpha\1(x\tl y\tl z)&=(x\tl y)\tl z,\\
\alpha\1(a\tb b\tb c)&=(a\t b)\tb c.
\end{align*}

Similarly there exists a unique morphism of square groups
$$
\beta=\beta_{A,B,C}: A\tl (B\tl C)\to A\tl B\tl C
$$
which is the canonical isomorphism $A\ee\t B\ee\t C\ee\to A\ee\t (B\ee\t
C\ee)$ on the ee-level and on the e-level satisfies
\begin{align*}
\beta(x\tl(y\tl z))&=x\tl y\tl z,\\
\beta(x\tl(b\tb c))&=(x|x)_H\tb b\tb c,\\
\beta(a\tb(b\t c))&=a\tb b\tb c.
\end{align*}
This is an isomorphism with inverse given by
\begin{align*}
\beta\1(x\tl y\tl z)&=x\tl (y\tl z),\\
\beta\1(a\tb b\tb c)&=a\tb (b\t c).
\end{align*}

From these isomorphisms one obtains that $\tl$ is indeed associative with
associativity isomorphisms given as in Theorem \ref{tlmonoidal}. The last statement of
Theorem \ref{tlmonoidal} follows from the fact that one has
$$
x\tl (b\tb c)=x\tl P(b\t c)=(x|x)_H\tb (b\t c)
$$
and this corresponds to $((x|x)_H\t b)\tb c$.

\subsection{The pentagon axiom}
The well-known pentagon axiom for monoidal categories claims that two natural
ways from $((A\tl B)\tl C)\tl D$ to $A\tl (B\tl (C\tl D))$
are equal. That is the following diagram commutes.
$$
\xymatrix@!C=4em{
&&(A\tl B)\tl(C\tl D)\ar[rrd]^{\alpha_{A,B,C\tl D}}\\
((A\tl B)\tl C)\tl D\ar[rru]^{\alpha_{A\tl B,C,D}}\ar[rd]_{\alpha_{A,B,C}\tl D}&&&&A\tl(B\tl(C\tl D))\\
&(A\tl(B\tl C))\tl D\ar[rr]_{\alpha_{A,B\tl C,D}}&&A\tl((B\tl C)\tl
D)\ar[ur]_{A\tl\alpha_{B,C,D}} }
$$
To show
this statement in our circumstances let us observe that $(((A\tl B)\tl C)\tl
D)\e$ is generated by elements of the form $((x\tl y)\tl z)\tl w$, $((a\tb
b)\tl z)\tl w$, $(((a\t b)\tb c)\tl w$ and $(((a\t b)\t c)\tb d$. Then the
statement is clear for elements of the form $((x\tl y)\tl z)\tl w$ and
$(((a\t b)\t c)\tb d$. It is also straightforward to check that both ways
carry $(((a\t b)\tb c)\tl w$ to $a\tb (b\t (c \t \Delta (w)))$. It remains to
consider $((a\tb b)\tl z)\tl w$. By one way it goes to $a\tb (b\t
(\Delta(z)\t \Delta(w))$, while by the second way it goes to $a\tb (b\t
\Delta(z\tl w))$. Thus the pentagon axiom follows from the following Lemma.

\begin{Le}\label{tldelta}
For the homomorphism
$$
\Delta_{M\tl N}:\Cok(P_{M\tl N})=\Cok(P_M)\t \Cok(P_N)\to (M\tl N)\ee=M\ee\t N\ee
$$
one has
$$
\Delta_{M\tl N}=\Delta_M\t\Delta_{N}.
$$
\end{Le}
\begin{proof}
Since $a\tb b=P(a\t b)$ we have $\Delta(a\tb b)=0$. On the other hand
\begin{align*}
\Delta(x\tl y)=&(x\tl y|x\tl y)_H-H(x\tl y)+TH(x\tl y)\\
=&(x|x)_H\t (y|y)_H-(x|x)_H\t H(y)-H(x)\t \Delta(y)\\
&+T((x|x)_H\t H(y)+H(x)\t\Delta(y))\\
=&(x|x)_H\t (y|y)_H-(x|x)_H\t H(y)-H(x)\t ((y|y)_H-H(y)+TH(y))\\
&-T((x|x)_H)\t TH(y)-TH(x)\t T\Delta(y)\\
=&(x|x)_H\t (y|y)_H-(x|x)_H\t H(y)-H(x)\t (y|y)_H\\
&+H(x)\t H(y)-H(x)\t TH(y)+(x|x)_H\t TH(y)+TH(x)\t \Delta(y)\\
=&(x|x)_H\t (y|y)_H-(x|x)_H\t H(y)-H(x)\t (y|y)_H\\
&+H(x)\t H(y)-H(x)\t TH(y)\\
&+(x|x)_H\t TH(y)+TH(x)\t ((y|y)_H-H(y)+TH(y))
\end{align*}
and
\begin{align*}
\Delta(x)\t \Delta(y)=&((x|x)_H-H(x)+TH(x))\t ((y|y)_H-H(y)+TH(y))\\
=&(x|x)_H\t (y|y)_H-H(x)\t (y|y)_H+TH(x)\t(y|y)_H-(x|x)_H\t H(y)\\
&+H(x)\t H(y)-TH(x)\t H(y)+(x|x)_H\t TH(y)\\
&-H(x)\t TH(y)+TH(x)\t TH(y)
\end{align*}
which coincides with $\Delta(x\tl y)$.
\end{proof}

\subsection{Unit object for $\tl$}

We start to check that $\Znil$ has the unit object property, where
$(\Znil)\e=(\Znil)\ee=\ZZ$, $P=0$ and $H(n)=\binom n2$. For any square
group $A$ we define
$$
\iota(A)\e:(\Znil\tl A)\e\to A\e
$$
by
$$
\iota(A)\e(n\tl x)=nx+\binom n2PH(x)
$$
for $n\in\ZZ$, $x\in A\e$ and
$$
\iota(A)\e(n\tb a)=nP(a).
$$

Let us show that $\iota(A)\e$ respects all relations of Definition
\ref{tldef}. The relation \eqref{bice} of Definition \ref{tldef} is clear,
because $P$ is a homomorphism with values in the center of $A\e$. We have
$$
\iota(A)\e(n\tl (x_1+x_2))=n(x_1+x_2)+\binom n2 PH(x_1+x_2).
$$
Since $n(x_1+x_2)=nx_1+nx_2-\binom n2 [x_1,x_2]$ and $PH(x_1+x_2)=
PH(x_1)+PH(x_2)+[x_1,x_2]$, the relation \eqref{ledi} of Definition
\ref{tldef} follows. We also have
\begin{multline*}
\iota(A)\e((n_1+n_2)\tl x)=(n_1+n_2)x+\binom {n_1+n_2}2PH(x)\\
=n_1x+n_2x+\binom {n_1}2PH(x)+\binom {n_2}2PH(x)+n_1n_2PH(x)\\
=\iota(A)\e(n_1\tl x)+ \iota(A)\e(n_2\tl x)+\iota(A)\e(n_1n_2\tb H(x))
\end{multline*}
and the relation \eqref{ridi} of Definition \ref{tldef} follows, because for
$H(n)=\binom n2$ one has $(n_1,n_2)_H=n_1n_2$. Similarly, we have
\begin{align*}
\iota(A)\e(n\tl P(a))
&=nP(a)+\binom n2PHP(a)\\
&=nP(a)+\binom n2 2P(a)\\
&=n^2P(a)\\
&=\iota(A)\e(n^2\tb a)
\end{align*}
and the condition \eqref{rip} of Definition \ref{tldef} follows. Since $P=0$
for $\Znil$ we have
$$
\iota(A)\e(P(n)\tl x)=0=nP\Delta(x)=\iota(A)\e(n\tb \Delta(x))
$$
and the condition \eqref{lep} of Definition \ref{tldef} follows. Finally, one
has
$$
\iota(A)\e(T(n)\tb T(a))= -nPT(a)=-nP(a)=-\iota(A)\e(n\tb a)
$$
and the relation \eqref{TT} of Definition \ref{tldef} follows.

We now define
$$
\iota(A)\ee:(\Znil\tl A)\ee\to A\ee
$$
by
$$
\iota(A)\ee(n\ox a)=na.
$$
We claim that $\iota(A)=(\iota(A)\e,\iota(A)\ee)$ defines the natural
morphism of square groups $\iota(A): \Znil\tl A\to A$. Indeed, we have
$$
\iota(A)\e(P(a\t a))=\iota(A)\e(n\tb a)=nP(a)=P(na)=P(n\t a)
$$
and compatibility with $P$ follows. We have also
\begin{align*}
\iota(A)\ee H(n\tl x)
&=\iota(A)\ee (n^2\t H(x)+\binom n2 \t \Delta(x))\\
&=n^2H(x)+\binom n2 \Delta(x)\\
&=n^2H(x)+\binom n2 HPH(x)-\binom n2 2Hx+\binom n2 \cref Hxx\\
&=nH(x)+\binom n2 \cref Hxx+\binom n2 HPHx\\
&=H(nx+\binom n2 PHx)\\
&=H\iota(A)\e (n\tl x)
\end{align*}
and the claim follows.

Next we show that $\iota(A):\Znil\tl A\to A$
is an isomorphism. The inverse is given by
$$
\iota\1\e(A)(x)=1\tl x, \ \ \iota\1\ee(A)(a)=1\t a.
$$
Since
$$
P\iota\1\ee(A)(a)=P(1\t a)=1\tb a=1\tl Pa=\iota\1\e(A)(Pa)
$$
and
$$
H\iota\1\e(A)(x)=H(1\tl x)=1\t Hx= \iota\1\ee(A)(Hx)
$$
it follows that $\iota\1(A):A\to\Znil\tl A$ is indeed a morphism of square
groups. Since
\begin{align*}
&\iota\ee(A)\iota\1\ee(A)(a)=\iota\ee(A)(1\t a)=a,\\
&\iota\e(A)\iota\1\e(A)(x)=\iota\e(A)(1\tl x)=x,\\
&\iota\1\ee(A)\iota\ee(A)(n\t x)=1\t nx=n\t x
\intertext{and}
&\iota\1\e(A)\iota\e(A)(n\tl x)=1\tl nx+1\tl \binom n2 Phx=n\tl x
\end{align*}
we see that $\iota(A)$ is really an isomorphism. In a similar way we will see
that
$$
\ka(A):A\tl\Znil\to A
$$
is an isomorphism of square groups, where $A$ is an arbitrary square group and
$$
\ka(A)\e(x\tl n)=nx, \ \ \ka(A)\e(a\tb n)=nP(a), \ \ \ka(A)\ee(a\t n)=na.
$$

Finally, we have to check that
$$
\xymatrix{(A\tl\Znil)\tl B\ar[r] \ar[dr]_{\ka\tl \Id}& A\tl(\Znil\tl B)
\ar[d]^{\Id\tl \iota}\\
&A\tl B}
$$
is a commutative diagram, where the top map is given by the associativity
isomorphisms. Generators of $((A\tl\Znil)\tl B)\e$ are of the form $(x\tl
n)\tl y$, $(a\t n)\tl y$ and $(a\t n)\tb b$. Commutativity of the diagram is
obvious for $(a\t n)\tb b$. On the other hand for $(a\t n)\tl y$ it means
$nP(a)\tl y=a\tb n\Delta(a)$, which follows from (5) of Definition
\ref{tldef}. Finally, commutativity of the diagram for element $(x\tl n)\tl
y$ means $(nx)\tl y=x\tl (nx+\binom n2 PHy)$, which can be checked as
follows. By the identity \eqref{fminusis} and (2) of Definition \ref{tldef}
we have
$$
(nx)\tl y=n(x\tl y)+\binom n2 (x|x)_H\tb H(y)
$$
while $x\tl (nx+\binom n2
PHy)=n(x\tl y)+\binom n2 x\tl PHy=n(x\tl y)+\binom n2 (x|x)_H\tb
H(y)$. Here we used (2) and (4) of Definition \ref{tldef}. Now the proof that
 $\Square$ is a monoidal category is complete.

\subsection{Symmetry property of $\tl$}
Next we prove that the tensor product $\tl$ is symmetric monoidal. To this
end we define
$$
\tau(A,B)\e(x\tl y)=y\tl x-H(y)\tb TH(x)
$$
for $x\in A\e$, $y\in B\e$ and
$$
\tau(A,B)\ee(a\ox b)=b\ox a
$$
for $a\in A\ee$, $b\in B\ee$, which then also necessarily determines
$$
\tau(A,B)\e(a\tb b)=b\tb a
$$
and makes compatibility with $P$ clear. Compatibility with $H$ means
$$
H(y\tl x-H(y)\tb TH(x))=H(y)\ox(x|x)_H+\Delta(y)\ox H(x).
$$
Indeed we have
\begin{align*}
H(y\tl x&-H(y)\tb TH(x))\\
=&H(y\tl x)-H(H(y)\tb TH(x))+(H(y)\tb TH(x)|H(y)\tb TH(x))_H\\
=&(y|y)_H\ox H(x)+H(y)\ox\Delta(x)\\
&+HPH(y)\ox HPTH(x)-HPH(y)\ox TH(x)-H(y)\ox HPTH(x)\\
&+(P(H(y)\ox TH(x))|P(H(y)\ox TH(x)))_H\\
=&(\Delta(y)+H(y)-TH(y))\ox H(x)+H(y)\ox((x|x)_H-H(x)+TH(x))\\
&+HPH(y)\ox HPH(x)-HPH(y)\ox TH(x)-H(y)\ox HPH(x)\\
=&H(y)\ox(x|x)_H+\Delta(y)\ox H(x)\\
&-TH(y)\ox H(x)+H(y)\ox TH(x)\\
&+HPH(y)\ox HPH(x)-HPH(y)\ox TH(x)-H(y)\ox HPH(x),
\end{align*}
so compatibility of $\tau$ with $H$ amounts to showing that the sum
$$
-TH(y)\ox H(x)+H(y)\ox TH(x)+HPH(y)\ox HPH(x)-HPH(y)\ox TH(x)
-H(y)\ox HPH(x)
$$
is zero. Substituting here $HP=1+T$ gives
\begin{multline*}
-TH(y)\ox H(x)+H(y)\ox TH(x)+H(y)\ox H(x)+H(y)\ox TH(x)+TH(y)\ox H(x)\\
+TH(y)\ox TH(x)-H(y)\ox TH(x)-TH(y)\ox TH(x)-H(y)\ox H(x)-H(y)\ox TH(x),
\end{multline*}
which is indeed zero.

Naturality of $\tau$ is straightforward; to check $\tau(B,A)\tau(A,B)=\Id$,
the only nontrivial part is to look at
\begin{multline*}
\tau(B,A)\e\tau(A,B)\e(x\tl y)
=\tau(B,A)\e(y\tl x-H(y)\tb TH(x))\\
=x\tl y-H(x)\tb TH(y)-TH(x)\tb H(y).
\end{multline*}
But
\begin{align*}
-H(x)\tb TH(y)-TH(x)\tb H(y)
&=-P(H(x)\ox TH(y))-P(TH(x)\ox H(y))\\
&=-P(H(x)\ox TH(y))-P(TH(x)\ox TTH(y))\\
&=-P(H(x)\ox TH(y))+PT(H(x)\ox TH(y))\\
&=0.
\end{align*}
Finally we have to check that two hexagons commute. The first case amounts to
checking that two ways from $(A\tl B)\tl C$ to $B\tl(C\tl A)$ are the same.
That is the following diagram commutes.
$$
\xymatrix{&(B\tl A)\tl C\ar[r]& B\tl (A\tl C)\ar[rd]&\\
(A\tl B)\tl C\ar[ru]\ar[rd]&&&B\tl (C\tl A)\\
&A\tl (B\tl C)\ar[r]& (B\tl C)\tl A\ar[ru]&
}
$$
Again this is trivial for elements of the form $(a\t b)\tb c$. It is
straightforward to check that both ways take the element $(a\tb b)\tl z$ to
$b\tb (\Delta(z)\t a)$. For the element $(x\tl y)\tl z$ it amounts to showing
that $y\tl (z\tl x)-(y|y)_H\tb (Hz\t THx)-Hy\tb (\Delta(z)\t THx)= y\tl (z\tl
x)-y\tl (Hz\tb THx)-Hy\tb (\Delta(z)\t THx)$ which is obvious since of $y\tl
(Hz\tb THx)=y\tl P(Hz\t THx)=(y|y)_H\tb (Hz\t THx)$. The second hexagon axiom
amounts to showing that two ways from $A\tl (B\tl C)$ to $(C\tl A)\tl B$ are
the same. That is the following diagram commutes.
$$
\xymatrix{&A\tl(C\tl B)\ar[r]& (A\tl C)\tl B\ar[rd]&\\
A\tl(B\tl C)\ar[ru]\ar[rd]&&& (C\tl A)\tl B\\
&(A\tl B)\tl C\ar[r]& C\tl (A\tl B)\ar[ru]&
}
$$
This is trivial for $a\tb (b\t c)$. One checks that both ways
take the element $x\tl (b\tb c)$ to $(c\t (x|x)_H)\tb b$. Finally the element
$x\tl(y\tl z)$ goes to
$$
c\tl (a\tl b)-Hc\tb TH(a\tl b)=c\tl (a\tl b)-Hc\tb ((a|a)_H\t TH(b))-Hc\tb (THa\t
\Delta(b))
$$
in $C\tl (A\tl B)$ and to
$$
(c\tl a)\tl b- (Hc\t (a|a)_H)\tb TH(b)-(Hc\t THa)\tb\Delta(b))
$$
in $(C\tl A)\tl B$. By the second way the element $x\tl (y\tl z)$ goes to
$$
(c\tl a)\tl b-(Hc\tb THa)\tl b -(Hc\t (a|a)_H)\tb TH(b)
$$
But these elements are the same, because
$$
(Hc\tb THa)\tl b=P(Hc\t THa)\tl b= (Hc\t THa)\tb \Delta(b).
$$

\subsection{The action of the monoidal category $(\SG, \tl)$ on the category
of pre-square groups}\label{tlaction}

In Section \ref{tlpsg} we defined the bifunctor
$$
\tl:\PSG\x\Square\to\PSG
$$
which yields a right action of the monoidal category $(\SG,\tl)$ on the category
$\PSG$, meaning that for any $M\in \PSG$ and $N_1,N_2\in \SG$ there are coherent isomorphisms
\begin{align*}
(M\tl N_1)\tl N_2&\cong M\tl(N_1\tl N_2),\\
M\tl\Znil&\cong M.
\end{align*}
The proof of this fact is quite similar to the proof of
associativity and unit properties of
$(\Square,\tl)$, therefore we omit it.

Existence of the bifunctor $
\tl:{\PSG}\x {\Square}\to {\PSG}$ is crucial in the following
Lemma. First
observe that if $M$, $M'$ and $N$ are square groups and $f:\wp(M)\to \wp(M')$
is a morphism
of underlying pre-square groups, then $f\tl\Id$ defines a morphism of pre-square groups
$$
\wp(M\tl N)\to \wp(M\tl N').
$$
We can take $M'=M$ and $f=n^*$ (see Lemma \ref{intpsg}).
The following result shows that the action of integers respects the tensor product.

\begin{Le}\label{tln} Let $M$ and $N$ be square groups. Then for any integer $n$ one has
$$
(n^*)_M\tl \Id_N=(n^*)_{M\tl N}.
$$
\end{Le}
\begin{proof}
The result is obvious on the ee-level, while on the e-level it can be checked as follows
\begin{multline*}
(n^*(x))\tl y=(nx+P(\binom n2 Hx))\tl y=(nx)\tl y+P(\binom n2 Hx)\tl y\\
=(nx)\tl y+(\binom n2 Hx)\tr \Delta(y)=n(x\tl y)+ \binom n2 ((x|x)_H\tr
H(y)+H(x)\tr \Delta(y))\\
=n(x\tl y)+\binom n2 H(x\tl y)
\end{multline*}
Here we used the fact that $(nx)\tl y=n(x\tl y)=\binom n2 (x|x)_H\tb H(y)$,
which follows from identity \eqref{fminusis} in section \ref{4444}.
\end{proof}

%\section{Properties of the tensor product}

\section{The tensor product of abelian square groups and quadratic $\ZZ$-modules}

We now describe the tensor product of those square groups which are abelian or are
quadratic $\ZZ$-modules. In this case the tensor product has a particularly simple
form.

\begin{Pro}\label{qzmodsg}
Let $A$ be an abelian square group and $B$ be a square group, then $A\tl B$
is an abelian square group. Moreover the abelian group $(A\tl B)\e$ is given
by the following pushout diagram of abelian groups.
$$
\xymatrix{
A\ee\t \Cok(P^B)\ar[r]^-{\Id\t \Delta}\ar[d]^{P\t \Id}&
A\ee\t B\ee/{(\Id+T\t T)}\ar[d] \\
A\e\t \Cok(P^B) \ar[r]& (A\tl B)\e
}
$$
If additionally $A$ is a quadratic $\ZZ$-module then $A\tl B$ is a
quadratic $\ZZ$-module as well.
\end{Pro}

\begin{proof}
Since the cross-effect of $H$ vanishes on $A$ it follows that $\rho=0$, where
$\rho$ is the same as in Definition \ref{tldef}. Hence $H^{A\tl B}$ is a
homomorphism and therefore $A\tl B\in \Ab(\Square)$. As a consequence $(A\tl
B)\e$ is an abelian group. Moreover for each $x\in A\e$ one has $x\tl P(b)=0$
for all $b\in B\ee$, therefore the function $x\tl -$ factors through
$\Cok(P_A)$. Finally, the relations (5) and (6) of Definition \ref{tldef}
show that $(A\tl B)\e$ indeed fits into the above pushout diagram. If
additionally $\Delta^A=0$, then $\Delta^{A\tl B}=0$ (see Lemma
\ref{tldelta}), hence $A\tl B$ is a quadratic $\ZZ$-module.
\end{proof}

%\begin{Pro}\label{qzmodsg}
%Let $A$ be a quadratic $\ZZ$-module and $B$ be a square group. Then
%$(A\tl B)\e$ is the abelian group given by in the following
%pushout diagram of abelian groups.
%$$\xymatrix{A\ee\t \Cok(P^B)\ar[r]^-{\Id\t \Delta}\ar[d]^{P\t \Id}&
%A\ee\t B\ee/{(\Id+T\t T)}\ar[d] \\
%A\e\t \Cok(P^B) \ar[r]& (A\tl B)\e}
%$$
%\end{Pro}

\begin{Pro}\label{qzmods}
Let $A$ be an abelian square group and $B$ be a quadratic $\ZZ$-module. Then
$$
A\tl B\cong E(A\ee\t B\ee,T\t T)\oplus \Cok(P^A)\t \Cok(P^B)
$$
is a quadratic $\ZZ$-module. Here the abelian group
 $\Cok(P^A)\t \Cok(P^B)$ is considered as a
square group by putting $0$ on the $\mathrm{ee}$-level and $E(-,-)$ is defined in Section \ref{qvazimodo}.
\end{Pro}

\begin{proof}
This is clear because $\Delta^B=0$.
\end{proof}

\begin{Le}\label{ttl} Let $A$ be an abelian group and $M$ be a square group.
Then
$$
f:(A\t M\ee)^\t\to A^\t \tl M
$$
is an isomorphism. Here
$$
f\e(a\t c)=(a,0)\tb c
$$
and
$$
f\ee(a\t c, b\t d)=(a,0)\t c - (0,b)\t T(d),
$$
with $a,b\in A$ and $c,d\in M\ee$.
\end{Le}

\begin{proof} The statement is obvious on the ee-level. To see it
on the e-level one observes that $P$ is a split epimorphism on $A^\t$ and
$\Ker(P)$ consists of elements $(a,-a)\in A^\t\ee$, $a\in A$. By Proposition
\ref{qzmodsg} one can check that $(A^\t\tl M)\e$ is the quotient of $(A\oplus
A)\t M\ee$ by the relations
\begin{align*}
(a,-a)\t \Delta(x)&=0, &\textrm{$a\in A$, $x\in M\e$,}\\
(a,0)\t c+(0,a)\t T(c)&=0, &c\in M\ee.
\end{align*}
Since $T\Delta +\Delta=0$ we see that the first relation follows from the
second one. Now the result is clear.

\end{proof}
%%%%%%%%%%%%%%%%%%%%%%%%%%%%%%%%%%%%%%%%%%%%%%%%
\section{The tensor product $V(S)\tl M$}
%%%%%%%%%%%%%%%%%%%%%%%%%%%%%%%%%%%%%%%%%%%%%%%%%
We consider the case when one of the factors in a tensor product is the free
square group $V(S)$ (see Section \ref{Vconst}). We take $S=\{1,2,\cdots,
n\}$. In this case we use the notation $V(n)$ instead of $V(S)$. If $n=1$
then $V(1)\cong \ZZ^Q$ and we have the following explicit result.

\begin{Le} For any square group $A$ one has an isomorphism
$$
\ZZ^Q\tl A \cong (A\ee\x A\e\xto{H^{\ZZ^Q\tl A}} A\ee\oplus
A\ee\oplus A\ee\xto{P^{\ZZ^Q\tl A}} A\ee\x A\e)
$$
where
$$
P^{\ZZ^Q\tl A}(a,b,c)= (a+b-Tb,Pa)
$$
and
$$
H^{\ZZ^Q\tl A}(a,x)=(a+T(a)+\Delta(x), -Ta, H(x)).
$$
\end{Le}

\begin{proof}
We just indicate the explicit isomorphism $\psi:A\e\x A\ee\to (\ZZ^Q\tl A)_c$
and its inverse:
\begin{align*}
\alpha\e(x,a)&=(1,0)\tl x+(1,0,0)\tb a,\\
\alpha\1\e((m,n)\tl x+(k_1,k_2,k_3)\tb a)&=(mx+\binom m2 PHx+k_1Pa,
n\Delta(x)+ k_1a+k_2a-k_2Ta).
\end{align*}
\end{proof}

For arbitrary $n$ and a square group $M$ let $G_n(M)$ be the group, which is
$(M\e)^n\x (M\ee)^{\binom n2}$ as a set with $(x_k,a_{ij})$ as a generic
element. Here $1\le k\le n$ and $1\le i<j\le n$. The group structure is given
by
$$(x_k,a_{ij})+(y_k,b_{ij})=(x_k+y_k,a_{ij}+b_{ij}+\cref H{x_j}{x_i}).$$
\begin{Pro}\label{vstl}
For any integer $n\ge1$ and any square group $M$ one has isomorphisms of
groups
\begin{align*}
(V(n)\tl M)\ee&\cong (M\ee)^{n^2+2n},\\
(V(n)\tl M)\e&\cong G_n(M)\x(M\ee)^n.
\end{align*}
\end{Pro}

\begin{proof} The isomorphism is obvious on the
$\mathrm{ee}$-level. We define
$$
\alpha\e:G_n(M)\x (M\ee)^n\to (V(n)\tl M)\e
$$
as follows. If a generic element of the group $G_n(M)\x (M\ee)^n$ is
$(x_k,a_{ij},b_l)$, where $1\le k,l\le n$, $1\le i<j\le n$ and $x_k\in M\e$,
$a_{ij},b_l\in M\ee$, then we put
$$
\alpha\e (x_k,a_{ij},b_l)=1\tl x_1+2\tl x_2+\cdots +n\tl x_n+\sum_{i<j}(j,i)\tb
a_{ij}+\sum_{l}Hl\tb b_l.
$$
One easily checks that $\alpha$ is a homomorphism. Actually it is an
isomorphism with inverse $\alpha\1\e$ which is uniquely determined by
\begin{multline*}
\alpha\1\e(\sum_iPHi\tl x_i+\sum_lHl\tb d_l+\sum_{ij}(i,j)\tb
c_{ij}+\sum_kHPHk\tb
b_k)=\\
(PHc_{kk},c_{ji}+Tc_{ij},d_l+b_l-Tb_l+\Delta(x_l)),
\end{multline*}
and
\begin{multline*}
\alpha\1\e(m_11+\cdots +m_nn+\sum_{i<j}m_{ij}[i,j])\tl x=\\
(m_kx+\binom {m_k}2 PHx,n_im_jHx+m_{ij}(Hx-THx),0).
\end{multline*}
\end{proof}

\section{Right exactness}\label{rexacttl}

We prove Proposition \ref{tlright}. First we check that for
any square group $A$ the tensor product functor $A\tl-:\SG\to\SG$ preserves
reflexive coequalizers. Recall that a functor $R:\Square \to \Square$
preserves reflexive coequalizers if for any simplicial object $B_*$ in
$\Square$ the canonical morphism $\pi_0(R(B_*))\to R(\pi_0(B_*))$ is an
isomorphism.
%Here $\pi_*(B_*)\in \Square$ is defined by level-wise homotopy
%groups:
%$$
%\pi_*(B_*)\e:=\pi_*({B_*}\e) \ \ \pi_*(B_*)\ee:=\pi_*({B_*}\ee)
%$$
%while $P$ and $H$ on $\pi_*(B_*)$ are induced from $B_*$.
Observe that $\pi_0(B_*)$ is a coequalizer of two parallel arrows
$d_1,d_0:B_1\to B_0$. We put $B=\pi_0(B_*)$. By universality property of
coequalizers we have a canonical map $\pi_0(A\tl B_*)\to A\tl B$. This is an
isomorphism, because it has an inverse, which is defined as follows. Take a
generator $x\tl y$ with $x\in A\e$ and $y\in B\e$. We choose an element
$z\in{B_0}\e$ in the class $y\in\pi_0({B_0}\e)$. The class of $x\tl\hat{y}\in
A\tl B_0$ in the quotient $\pi_0(A\tl B_*)$ is independent of the choice.
Indeed, if $z'$ is also in the class $y$, then there exists $w\in{B_1}\e$
such that $d_0w=z$ and $d_1w=z'$ and therefore $x\tl z$ and $x\tl z'$ define
the same element in $\pi_0(A\tl B_*).$ Based on this fact one easily checks
that this assignment respects all relations for $\tl$ and indeed defines a
morphism $A\tl B\to \pi_0(A\tl B_*)$ and hence $A\tl-:\SG\to \SG$ preserves
reflexive coequalizers.

Next we prove that the tensor product $A\tl-:\SG\to \SG$ preserves finite
products. Let $B,C\in \SG$. Then $B\x C$ in $\SG$ is constructed degreewise,
i.~e. $(B\x C)\e=B\e\x C\e$, $(B\x C)\ee=B\ee\x C\ee$, $P(b,c)= (Pb,Pc)$ and
$H(y,z)=(Hy,Hz)$. Here $y\in B\e,z\in C\e$ and $b\in B\ee, c\in C\ee$. The
projection $p_1:B\x C\to B$ has the canonical section $i_1:B\to B\x C$ given
by ${i_1}\e(x)=(x,0)$ and ${i_1}\ee(b)=(b,0)$. This section yields a morphism
of square groups $\Id_A\tl i_1:A\tl B\to A\tl (B\x C)$. Similarly, one gets
the morphism $\Id_A\tl i_2:A\tl C\to A\tl (B\x C)$. The identity
\eqref{tlcom} in Section \ref{tlpsg} shows that the images of $\Id_A\tl i_1$ and $\Id_A\tl i_2$
commute and therefore they yield the canonical morphism $$i_*:(A\tl B)\x
(A\tl C)\to A\tl (B\x C)$$ which obviously is left inverse to the canonical
morphism $A\tl (B\x C)\to (A\tl B)\x (A\tl C)$ induced by projections. It is
clear that both morphisms are isomorphisms on the ee-level. Thus we need only
to show that the map $i_*$ is surjective on the e-level. This follows
immediately from the fact that $x\tl (y,z)=x\tl i_1(y)+x\tl i_2(z)$ and $a\tl
(b,c)=a\tl i_1(b)+a\tl i_2(c)$.

We now prove that for any short exact sequence of square groups
$$
0\to B_1\xto{\mu} B\xto{\sigma} B_2\to 0
$$
the induced sequence
$$
A\tl B_1\to A\tl B\to A\tl B_2\to 0
$$
is also exact. To this end one observes that by Lemma \ref{centralizator} it
suffices to consider the case, when the extension is central. Consider the
following diagram in $\SG$
$$\xymatrix{
B_1\x B\ar@<.5ex>[r]^-f\ar@<-.5ex>[r]_-g&B}
$$
where $g$ is the projection on the second factor, while $f\e(b_1,b)=\mu\e(b_1)+b$
and $f\ee(y_1,y)=\mu\ee(y)+y$, $y\in B\e,y_1\in B_{1e}$, $b_1\in B_{1ee},
b\in B\ee$. Since $B_1$ is central in $B$, it follows that $f$ is a morphism in
$\Square$ and in fact $B_2$ is isomorphic to the coequalizer of this diagram.
Since $A\tl -$ preserves products and reflexive coequalizers we see that
$$\xymatrix{
A\tl B_1\x A\tl B\ar@<.5ex>[r]^-{ f_*}\ar@<-.5ex>[r]_-{ g_*}&A\tl B\ar[r]& A\tl B_2} $$
is a coequalizer. From this follows that the functor $A\tl (-)$ is right exact.

Finally, assume that $A$ is a projective square group. We have to show that
the functor $A\tl (-)$ is exact. We can assume that $A=V(S)\x (\ZZ[T])^\t$
for some $S\in \Sets$ and $T\in \Sets$, because any projective is a retract
of a free one. Since the tensor product commutes with filtered colimits we
can assume that $S$ and $T$ are finite sets. Since the tensor product
commutes with finite products the result follows from Lemma \ref{ttl} and
Proposition \ref{vstl}.

\section{$\Znil[-]$ as a monoidal functor}
The next result generalizes the well-known fact that the free abelian group functor
$$\ZZ[-]:(\Sets,\x)\to (\Ab, \t)
$$
is a symmetric monoidal functor.
\begin{Pro}\label{tlsets}
For any sets $S$ and $S'$, one has a natural isomorphism of square
groups
$$
\delta:\Znil[S\x S']\cong\Znil[S]\tl\Znil[S']
$$
which on the $\mathrm{ee}$-level is the canonical isomorphism $\ZZ[S\x S'] \cong
\ZZ[S]\t \ZZ[S']$, given by $\delta\ee((s,s')\t (t,t'))=s\t t \t s'\t t'$.
On the $\mathrm e$-level it is given by
$$
\delta \e(s, s')=s\tl s', \ \ s,s'\in S.
$$
Thus the functor
$$
\Znil[-]:(\Sets,\x)\to (\SG, \tl)
$$
is symmetric monoidal.
\end{Pro}

\begin{proof}
Since $(\Znil[S\x S'])\e$ is a free nil$_2$-group on $S\x S'$, the
homomorphism $\delta \e$ is well-defined. Let us first check that the pair
$\delta=(\delta\e,\delta \ee)$ defines a morphism of square groups. We have
$$
P\delta \ee((s,s')\t (t,t'))=P(s\t t\t s'\t t')=(s\t t)\tb (s'\t t').
$$
Similarly we have
\begin{multline*}
\delta \e P((s,s')\t (t,t'))=\delta\e ([(t,t'),(s,s')])=
[\delta\e (t,t'),\delta\e (s,s')]=[t\tl t',s\tl s']\\
=(t|s)_H\t(t'|s')_H= s\t t\t s'\t t'.
\end{multline*}
Comparing these expressions we see that $\delta$ is compatible
with $P$. We have also
$$
\delta\ee ((s,s')|(t,t'))_H=\delta\ ((t,t')\t
(s,s'))=t\t s\t t'\t s'.
$$
On the other hand we have
$$
(\delta\e (s,s')\mid \delta\e(t,t'))_H=
(s\tl s'| t\tl t')_H=(s|t)_H\t (s'|t')_H=t\t s\t t'\t s'.
$$
Since
$$
\delta\ee H(s,s')=0
$$
and
$$
H(\delta\e((s,s'))=H(s\tl s')=(s|s')_H\t H(s')+H(s)\tl
\Delta(s')=0
$$
we conclude that $\delta$ is compatible with $H$ and therefore $\delta$ is,
in fact, a morphism of square groups. Since $\delta\ee$ is an isomorphism, it
suffices to show that $\delta\e$ is an isomorphism. We construct the inverse
map
$$
\eta:(\Znil[S]\tl\Znil[S'])\e
\to(\Znil[S\x S'])\e
$$
as follows. Take $s\in S$. Since $(\Znil[S'])\e$ is a
free nil$_2$-group on $S'$, there is a unique homomorphism
$$
f_{s}:(\Znil[S'])\e\to(\Znil[S\x S'])\e
$$
such that $f_{s}(s')=(s,s')$, $s'\in S'$.
In particular for any $y\in(\Znil[S'])\e$ the element $f_s(y)$
is well-defined. Since $(\Znil[S])\e$ is a free nil$_2$-group
we can extend $s\mapsto f_s(y)$ to a map $x\mapsto f_x(y)\in
(\Znil[S\x S'])\e$ in such a way that
$$
f_{x_1+x_2}(y)=f_{x_1}(y)+f_{x_2}(y)+g(x_1\t x_2\t H(y)),
$$
where $g:\ZZ[S]\t \ZZ[S]\t\ZZ[S']\t\ZZ[S']\to(\Znil[S\x S'])\e$ is given by
$g(s\t t\t s'\t t')= [(t,t'),(s,s')]$. Thus we have
$$
[f_{x_1}(y_1),f_{x_2}(y_2)]=-g(x_1,x_2,y_1,y_2).
$$
We claim that for all $x\in(\Znil[S])\e$ and $y_1,y_2\in(\Znil[S'])\e$ one has
$$
f_x(y_1+y_2)=f_x(y_1)+f_x(y_2).
$$
By our construction the claim holds if $x\in S$. Therefore it
suffices to show that the equation holds for $x=x_1+x_2$, provided
it holds for $x_1$ and $x_2$. To this end, we proceed as follows:
$f_{x_1+x_2}(y_1+y_2)=f_{x_1}(y_1)+f_{x_1}(y_2)+
f_{x_2}(y_1)+f_{x_2}(y_2)+g(x_1\t x_2\t H(y_1+y_2))=
f_x(y_1)+f_x(y_2)+[f_{x_1}(y_1),f_{x_2}(y_2)]+
g(x_1\t x_2\t y_2\t y_1) $ because $(y_1|y_1)_H=y_2\t y_1$ and last
two summands cancel; hence the claim.

Now we are ready to define the map
$\eta:(\Znil[S]\tl\Znil[S'])\e\to(\Znil[S\x
S'])\e$ by
\begin{align*}
\eta(x\tl y)&=f_x(y),\\
\eta((s\t t)\tb (s'\t t'))&=g(s\t t\t s'\t t').
\end{align*}
We have to check that $\eta$ respects the relations (1)-(6) of Definition
\ref{tldef}. The relations (1) and (6) are clear, (3) holds by our
construction and (2) we just checked. Let us check (4). Without loss of
generality we can assume that $b=s'\t t'$ with $s',t'\in S'$. First we
consider the case $x=s\in S$. Then one has $\eta(x\tl P(b))=\eta(s\tl
[t',s'])=[(s,t'),[(s,s')]=g(s\t s \t s'\t t')=\eta((s\t s)\tb s'\t
t')=\eta((x|x)_H\tb b)$. For general $x$ it suffices to show that $\eta$
respects the equality (4) for $x=x_1+x_2$ provided $\eta$ respects it for
$x_1$ and $x_2$. To check the last assertion one needs to show the equality
for the cross-effects of both sides of the equality in question. But this
is formal,
$$
\eta((x_1|x_2)_H\tb HP(b))=\eta((x_1|x_2)_H\tb (b+Tb))= \eta((x_1|x_2)_H\tb b+
(x_2|x_1)_H\tb b).
$$
Here we used the fact that $\eta$ respects the identity (6).
Now we check that $\eta$ respects the equation (5). Indeed, we can
assume that $y=s'\in S'$ and $a=s\t t$, $s,t\in S$. Then we have
$$
\eta(P(a)\tl y)=\eta([t,s]\tl s')=[(t,s'),(s,s')].
$$
On the other hand we have
$$
\eta(a\tb \Delta (b))=\eta(s\t s\t s'\t s')=[(t,s'),(s,s')]
$$
and the result follows.
\end{proof}

\section{Torsion product of square groups}\label{tortl}
It follows from Proposition \ref{tlright} that the functor $M\tl(-)$ respects
weak equivalences of simplicial square groups provided $M$ is projective. We
will exploit this fact in this section.

Thanks to Section \ref{ssg} for any square group $M$ one can take a cofibrant
replacement $M^c$ of $M$. This means that there is given a weak equivalence $M^c\to
M$ and $M^c$ is cofibrant in the model category structure introduced in Section
\ref{ssg}.

\begin{Le} For any square groups $M$ and $N$ there is an isomorphism of square
groups
$$\pi_i(M^c\tl N)\cong \pi_i(M\tl N^c), \ i\ge0$$
\end{Le}

\begin{proof}
According to \cite{qu} the square groups $\pi_i(M^c\tl N)$ (as well as
$\pi_i(M\tl N^c)$) do not depend on the cofibrant replacement. Moreover one
can assume that each component of $M^c$ and $N^c$ is a free square group
(\cite{qu}). One considers now the bisimplicial square group $M^c\tl N^c$.
Both on the e- and on the ee-level one has the Quillen spectral sequences for
double simplicial groups \cite{qss}. Since the tensor product with a free
object respects weak equivalences, both spectral sequences degenerate
yielding the isomorphism in question.
\end{proof}

We now put
$$
\ttl_i(M,N)=\pi_i(M^c\tl N)\cong \pi_i(M\tl N^c).
$$

\begin{Le}
For any square groups $M$ and $N$ one has natural isomorphisms
$\ttl_0(M,N)\cong M\tl N$ and $\ttl_1(M,N)\in \Ab(\Square)$ and for all
$i\ge2$ one has $\ttl_i(M,N)\in \Ab$.
\end{Le}

\begin{proof} The statement on $\ttl_0$ follows from right exactness of $\tl$.
As already mentioned, for any simplicial square group $X$ one has
$\pi_i(X)\in \Ab(\Square)$ for all $i\ge1$, hence the statement for $i=1$.
Finally let us observe that
$$
(\ttl_i(M,N))\ee=\mathsf{Tor^\ZZ}_i(M\ee,N\ee), \ \ i\ge0,
$$
hence $(\ttl_i(M,N))\ee=0$ for $i\ge2$ and the result follows.
\end{proof}

Thus we have bifunctors
\begin{align*}
\ttl_1:&\Square\x \Square \to\Ab(\Square)
\intertext{and}
\ttl_i:&\Square\x\Square \to \Ab, &i\ge2.
\end{align*}

\begin{Pro}\label{ttles}If
$$0\to N_1\to N\to N_2\to 0$$
is a short exact sequence of square groups then for any square group $M$ one
has a long exact sequence of square groups
\begin{multline*}
\cdots \to \ttl_2(M, N_2)\to \ttl_1(M,N_1)\to\ttl_1(M, N)\\
\to\ttl_1(M,N_2)\to M\tl N_1\to M\tl N\to M\tl N_2\to 0.
\end{multline*}
Furthermore $\ttl_i(M,N)=0$ provided $i\ge1$ and either $M$ or $N$ is
projective.
\end{Pro}

\begin{proof}
If $M$ is free, then one can take $M^c=M$. A projective object is a retracts
of a free square group and hence the statement on projective objects is
clear. For general $M$ we choose an $M^c$ with degreewise free square groups.
According to Proposition \ref{tlright} one has the short exact sequence of
simplicial objects
$$
0\to M^c\tl N_1\to M^c\tl N\to M^c\tl N_2\to 0
$$
yielding the long exact sequence for homotopy groups. Here we used the
well-known fact that epimorphisms of simplicial groups are Kan fibrations.
\end{proof}

\begin{Co}\label{margaliti} Let $A$ be an abelian group and $M$ be a square group. Then
$$\ttl_1(A^\t,M)\cong (\mathsf{Tot}^\ZZ_1(A,M\ee))^\t$$
and $\ttl_k(A^\t,M)=0$ for $k\geq 2$.
\end{Co}

\begin{proof} Take a short exact sequence
 $0\to F_1\to F_0\to A\to 0$ of abelian groups with free abelian group $F_0$. It
 yields a short exact sequence of square groups
 $$0\to F_1^\t\to F_0^\t \to A^\t\to 0.$$
 Since free square groups $F_0^\t$ and $F_1^\t$ are free square groups the
 result follows from Lemma
 \ref{ttl} and
 Proposition \ref{ttles}.
\end{proof}
Now we give the proof of Proposition \ref{tlco}.
\begin{proof}[Proof of Proposition \ref{tlco}]
We use the short exact sequence \eqref{sqcopr} in Section \ref{copr}.
Proposition \ref{ttles} shows that it yields the exact sequence of square
groups
\begin{multline*}
\ttl_1(M, A\vee B)\to \ttl_1(M, A)\x \ttl_1(M, B)\\ \to M\tl
(\Cok(P^A)\t
\Cok(P^B))^\t\xto{\mu} M\t (A\vee B)\to M\tl (A\x B)\to 0.
\end{multline*}
Let us observe that for any category $\C$ with coproducts and zero object and
for any functor $F:\C\to\Gr$ with $F(0)=0$ the natural homomorphism
$r:F(X\vee Y)\to F(X)\x F(Y)$ is always an epimorphism. Indeed, if $s_1\in
F(X)$ and $s_2\in F(Y)$ then $s=i_1(s_1)+i_2(s_2)\in F(X\vee Y)$ has the
property that $r(s)=(s_1,s_2)$. Here $i_1:X\to X\vee Y$ and $i_2:Y\to X\vee
Y$ are standard inclusions. Thus $\mu$ is a monomorphism and the result
follows from Proposition \ref{ttl} and the fact that $M\tl (-)$ preserves
finite products.
\end{proof}

\begin{Co} Let $M$ be a square group. Then the functor $\ttl_k(M,-)$ preserves finite products
for all $k\geq 0$. Moreover the natural map
$$\ttl_k(M,A\vee B)\to  \ttl_k(M,A)\x \ttl_k(M,B)$$
is an isomorphism provided $k\geq 2$. Here $A$ and  $B$ are square  groups. If $k=1$,
then one has a short exact sequence of square
groups
$$0\to (\mathsf{Tor^\ZZ_1}(M\ee,\Cok(P^A)\t \Cok(P^B)))^\t
\to \ttl_1(M,A\vee B)\to  \ttl_1(M,A)\x \ttl_1(M,B)\to 0$$
\end{Co}

\begin{proof} Since
$$\pi_*(M^c\tl (A\x B))=\pi_*(M^c\tl A)\x \pi_*(M^c\tl B)$$
we have
$$\ttl_*(M,A\x B)=
\ttl_*(M,A)\x \ttl_*(M,B)$$ and
the statement on finite products follows. For the
rest one apllies the long exact sequence for $\ttl$-functors to
the short exact sequence in
Proposition \ref{tlco} and
uses Corollary  \ref{margaliti}.
\end{proof}
\section{Quadratic rings}

A monoid in the monoidal category $(\SG,\tl)$ is termed \emph{quadratic
ring}. More explicitly, a quadratic ring structure on a square group $R$ is
given by a multiplicative monoid structure on $R\e$ and a ring structure on
$R\ee$. The multiplicative unit of $R\e$ is denoted by $1$. One requires that
these structures satisfy the following additional properties. First of all we
have
\begin{enumerate}
\item[(i)] $x(y+z)=xy+xz$,
\item[(ii)] $(x+y)z=xz+yz+P(\cref HyxH(z))$.
\end{enumerate}

Thus $\Cok(P)$ is a ring. Moreover the maps
\begin{align*}
-T&:R\ee\to R\ee\\
(-\mid -)_H&:\Cok(P_R)\t\Cok(P_R)\to R\ee
\end{align*}
are ring homomorphisms, in other words one has

\begin{enumerate}
\item[(iii)] $(x\mid y)_H (u\mid v)_H=(xu\mid yv)_H$,
\item[(iv)] $T(ab)+T(a)T(b)=0$.
\end{enumerate}

Let us observe that $T(abc)=T(a)T(b)T(c)$.

Furthermore the following equations hold
\begin{enumerate}
\item[(v)] $ P(a\Delta(x))=P(a)x,$
\item[(vi)] $P((x\mid x)_Ha)=xP(a),$
\item[(vii)] $H(xy)=(x\mid x)_HH(y)+H(x)\Delta(y).$
\end{enumerate}

It follows from Lemma \ref{tldelta} that
$$
\Delta:\Cok(P_R)\to R\ee
$$
is a ring homomorphism.

Let $\mathsf{QR}$ denote the category of quadratic rings. We have a full
embedding of categories
$$\mathsf{Rings}\subset\mathsf{QR}
$$
which identifies rings with quadratic rings $R$ satisfying $R\ee=0$. This
inclusion has a left adjoint given by $R\mapsto \Cok(P_R)$.

Let $R$ be a quadratic ring. A \emph{right $R$-quadratic module} is an object
$M\in \Square$ together with a right action of the monoid $R$ in the monoidal
category $(\Square, \tl)$. Equivalently, it is given by a square group $M$
together with a right $R\ee$-module structure on $M\ee$ and a right action of
the multiplicative monoid $R\e$ on $M\e$ such that the following holds (
$m,n\in M\e$, $c\in M\ee$)

\begin{enumerate}
\item[(i)] $m(x+y)=mx+my$,
\item[(ii)] $(m+n)x=mx +nx+P(\cref HnmH(x))$.
\item[(iii)] $(m\mid n)_H\cref Hxy=\cref H{mx}{ny}$,
\item[(iv)] $T(ca)+T(c)T(a)=0$,
\item[(v)] $ P(c\Delta(x))=P(c)x,$
\item[(vi)] $P(\cref Hmma)=mP(a),$
\item[(vii)] $H(mx)=\cref HmmH(x)+H(m)\Delta(x).$
\end{enumerate}
 The category of right $R$-quadratic modules is denoted by ${\sf QM}$-$R$. In a
 similar manner one can define the notion of a left $R$-quadratic module.

\subsection{Monoid quadratic rings}
Thanks to Proposition \ref{tlsets} the functor $\Znil[-]$ is a monoidal
functor from the monoidal category of sets to square groups $(\SG, \tl)$. It
follows that for any monoid $M$ one gets a quadratic ring $\Znil[M]$ and
thus the functor
$$
\mathsf{Monoids}\to\mathsf{QR}.
$$
 Now we observe that for a quadratic ring $R$
 the set of linear elements $\LL(R)$ is a multiplicative submonoid of
$R\e$. This follows directly from the formula $H(xy)=\cref
HxxH(y)+H(x)\Delta(y)$. Thus we obtain the functor
$$
\LL:\mathsf{QR}\to\mathsf{Monoids}
$$
The following result is a multiplicative version of Proposition
\ref{slineareleents}.

\begin{Pro}\label{lineareleents}
The functor $\Znil[-]:\mathsf{Monoids}\to\mathsf{QR}$ is left adjoint to the functor
$\LL$.
\end{Pro}

\begin{proof}
Let $M$ be a monoid and let $R$ be a quadratic ring. Given a morphism
$f:\Znil[M]\to R$ of quadratic rings, the corresponding map $f_0:M\to \LL(R)$
is just the composite $M\subset\Znil[M]\to R$ and hence is multiplicative.

Conversely we must show that any homomorphism $f_0:M\to\LL(R)\subset R\e$ of
monoids extends uniquely to a quadratic ring morphism $f:\Znil[M]\to R$.
In the proof of Proposition \ref{slineareleents} we already constructed
the morphism of square groups $f:\Znil[M]\to R$
extending $f_0$. It remains to show that $f\e$ and $f\ee$ are multiplicative. First
note that the induced map
$$
f_*:\ZZ[M]=\Cok(P_{\Znil[M]})\to\Cok(P_R)
$$
is multiplicative as it is the unique additive extension of the composite monoid
homomorphism
$$
M\to\LL(R)\to R\e\onto\Cok(P_R).
$$

It follows that $f\ee$ is multiplicative since by definition it factors into
the composition of multiplicative maps as follows.
$$
\ZZ[M]\ox\ZZ[M]=\Cok(P_{\Znil[M]})\ox\Cok(P_{\Znil[M]})\xto{f_*\ox f_*}
\Cok(P_R)\ox\Cok(P_R)\xto{(\_|\_)_H}
R\ee
$$
Finally the function $\Xi$ of two variables on $\Znil[M]\e$ given by
$$
\Xi(x,y)=f\e(xy)-f\e(x)f\e(y)
$$
is central in $R\e$. This function vanishes if both $x$ and $y$ are in $M$, so if we show that it
is biadditive it will follow that $f\e$ is multiplicative. Indeed
\begin{multline*}
\Xi(x,y+y')=f\e(xy+xy')-f\e(x)(f\e(y)+f\e(y'))=f\e(xy)+f\e(xy')-f\e(x)f\e(y')\\
=\Xi(xy)+\Xi(xy')
\end{multline*}
and
\begin{multline*}
\Xi(x+x',y)=f\e(xy+x'y+P((\bar x\ox\bar x')H(y)))-f\e(x+x')f\e(y)=\\
f\e(xy)+f\e(x'y)
+Pf\ee((\bar x\ox\bar x')H(y))-(f\e(x)f\e(y)+f\e(x')f\e(y)
+P((\bar{f\e(x)}\ox\bar{f\e(x')})Hf\e(y)))\\
=\Xi(x,y)+\Xi(x',y)+P(f\ee(\bar x\ox\bar x')f\ee
H(y))-(\bar{f\e(x)}\ox\bar{f\e(x')})Hf\e(y))=\Xi(x,y)+\Xi(x',y)
\end{multline*}
and we are done.
\end{proof}

\subsection{Commutative quadratic rings} Since
$\tl$ defines a symmetric monoidal category structure on $\SG$, one can talk
about commutative monoids in this monoidal category. We call them commutative
quadratic rings. Hence by definition a quadratic ring $R$ is commutative if
the following equations hold
\begin{align*}
ba&=ab, &a,b\in R\ee;\\
yx&=xy-P(H(x)TH(y)), &x,y\in R\e.
\end{align*}

In the following result we obtain a kind of \emph{cup$_1$-product.} We refer
to Remark \ref{comspectra} below for the homotopy theoretic meaning of the
groups involving it. The elements of the form $k^M(x)$ are the same as in
Corollary \ref{29}.

\begin{The}\label{psioperaciebi}
Let $R$ be a commutative quadratic ring, then the map
$R\e\to \Ker(P)$ of degree 4 given by $x\mapsto H(x)TH(x)\in \Ker(P)$ yields a well-defined quadratic map
$$
\psi:\Cok(P_R)\to\Ker(P:R\ee/(\Id-T)\to R\e)
$$
satisfying
$$
\psi(\bar{x}+\bar{y})=\psi(\bar{x})+\psi(\bar{y})-k^R(x)k^R(y)
$$
and
$$
\psi(\bar{x}\bar{y})=k^R(x^2)\psi(\bar{y})+\psi(\bar{x})k^R(y^2).
$$
\end{The}

\begin{proof}
We have
\begin{align*}
\psi(\bar x+\bar y)=&H(x+y)TH(x+y)\\
=&(H(x)+H(y)+(x|y)_H)(TH(x)+TH(y)-(y|x)_H)\\
=&H(x)TH(x)+H(y)TH(y)+H(x)TH(y)+H(y)TH(x)\\
&-H(x)(y|x)_H+(x|y)_HTH(x)-H(y)(y|x)_H+(x|y)_HTH(x)-(xy|yx)_H\\
=&\psi(\bar x)+\psi(\bar y)+H(x)TH(y)-T(TH(y)H(x))\\
&-H(x)(y|x)_H+T((y|x)_HH(x))-H(y)(y|x)_H+T((y|x)_HH(y))-(xy|yx)_H;
\end{align*}
using commutativity of $R\ee$, this is equal to
$$
\psi(\bar x)+\psi(\bar y)-(xy|yx)_H+(1-T)(H(x)TH(y)-H(x)(y|x)_H-H(y)(y|x)_H).
$$
On the other hand,
\begin{align*}
\Delta(x)\Delta(y)=&((x|x)_H-H(x)+TH(x))((y|y)_H-H(y)+TH(y))\\
=&(xy|xy)_H-(x|x)_HH(y)+(x|x)_HTH(y)-H(x)(y|y)_H+TH(x)(y|y)_H\\
&+H(x)H(y)+TH(x)TH(y)-H(x)TH(y)-TH(x)H(y)\\
=&(xy|xy)_H-(x|x)_HH(y)+T((x|x)_HH(y))\\
&-H(x)(y|y)_H+T(H(x)(y|y)_H)+H(x)H(y)-T(H(x)H(y))\\
&-H(x)TH(y)+T(H(x)TH(y))\\
=&(xy|xy)_H\\
&+(1-T)(-(x|x)_HH(y)-H(x)(y|y)_H+H(x)H(y)-H(x)TH(y)).
\end{align*}
Thus it remains to prove
$$
(xy|yx)_H=(xy|xy)_H;
$$
but by the hypothesis $yx$ is equal to $xy$ modulo image of $P$, and
$(\_|\_)_H$ vanishes if one of the operands is in the image of $P$.

For the second equality in the theorem one calculates
\begin{align*}
\psi(xy)=&H(xy)TH(xy)\\
=&((x|x)_HH(y)+H(x)\Delta(y))T((x|x)_HH(y)+H(x)\Delta(y))\\
=&((x|x)_HH(y)+H(x)\Delta(y))((x|x)_HTH(y)+TH(x)\Delta(y))\\
=&(x^2|x^2)_HH(y)TH(y)+H(x)TH(x)\Delta(y)^2\\
&+(x|x)_HH(y)TH(x)\Delta(y)+H(x)\Delta(y)(x|x)_HTH(y)\\
=&(x^2|x^2)_H\psi(y)+\psi(x)\Delta(y^2)\\
&+(x|x)_HH(y)TH(x)\Delta(y)-T((x|x)_HH(y)TH(x)\Delta(y)).
\end{align*}
\end{proof}

\section{Comparison of $\tl$, $\tr$ and $\odot$}\label{totr}
In this section we prove Proposition \ref{tl=tr=to}. Actually, we only show
that $\tl$ and $\odot$ are isomorphic. The rest is similar. Given the product
$A\tl B$ we define
$$
x\tr y=x\tl y-H(x)\tb TH(y).
$$
We then must check that the equalities \eqref{sridi} and \eqref{srip} hold in
$A\tl B$. Indeed we have
\begin{multline*}
(x_1+x_2)\tr y=(x_1+x_2)\tl y-H(x_1+x_2)\tb TH(y)\\
=x_1\tl y+x_2\tl y+(x_2|x_1)_H\tb H(y)-(H(x_1)+H(x_2)+(x_1|x_2)_H)\tb TH(y)\\
=x_1\tr y+x_2\tr y+(x_2|x_1)_H\tb H(y)-(x_1|x_2)_H\tb TH(y).
\end{multline*}
Since
m$$
(x_2|x_1)_H\tb H(y)=TT((x_2|x_1)_H)\tb TTH(y)=-T((x_2|x_1)_H)\tb TH(y)
=(x_1|x_2)_H\tb TH(y),
$$
We see that \eqref{sridi} indeed holds. For \eqref{srip} one considers
\begin{multline*}
x\tr P(b)=x\tl P(b)-H(x)\tb THP(b)=(x|x)_H\tb b-H(x)\tb HP(b)\\
=(x|x)_H\tb b-H(x)\tb b-H(x)\tb T(b)=(x|x)_H\tb b-H(x)\tb b-TTH(x)\tb T(b)\\
=(x|x)_H\tb b-H(x)\tb b+TH(x)\tb b=\Delta(x)\tb b.
\end{multline*}
Moreover we have
\begin{multline*}
H(x\tr y)=H(x\tl y-H(x)\tb TH(y))\\
=H(x\tl y)+H(-H(x)\tb TH(y))-(x\tl y|H(x)\tb TH(y))_H\\
=H(x\tl y)-H(H(x)\tb TH(y))+(H(x)\tb TH(y)|H(x)\tb TH(y))_H\\
=(x|x)_H\ox H(y)+H(x)\ox\Delta(y)-H(x)\ox TH(y)+TH(x)\ox H(y)\\
=((x|x)_H+TH(x))\ox H(y)+H(x)\ox(\Delta(y)-TH(y))\\
=(\Delta(x)+H(x))\ox H(y)+H(x)\ox((y|y)_H-H(y))
=\Delta(x)\ox H(y)+H(x)\ox(y|y)_H
\end{multline*}
%Finally
%$$
%\pi(x\tr y)=\pi(x\tl y-H(x)\tb TH(y))=\pi(x\tl y)
%$$
%as claimed.

We thus see that $A\tl B$ with $x\tr y$ defined as above has all properties
required of $A\odot B$ in Definition \ref{symdef}. Conversely, we must show that the
relations from Definition \ref{symdef} imply all the relations for $A\tl B$. We have
\begin{align*}
(x_1+x_2)\tl y=&(x_1+x_2)\tr y+H(x_1+x_2)\tb TH(y)\\
=&x_1\tr y+x_2\tr y+(H(x_1)+H(x_2)+(x_1|x_2)_H)\tb TH(y)\\
=&x_1\tl y-H(x_1)\tb TH(y)+x_2\tl y-H(x_2)\tb TH(y)\\
&+H(x_1)\tb TH(y)+H(x_2)\tb TH(y)+(x_1|x_2)_H\tb TH(y)\\
=&x_1\tl y+x_2\tl y+(x_1|x_2)_H\tb TH(y).
\end{align*}
Moreover
$$
(x_1|x_2)_H\tb TH(y)=-T(x_2|x_1)_H\tb TH(y)=(x_2|x_1)_H\tb H(y).
$$
Next, we check
\begin{multline*}
x\tl P(b)=x\tr P(b)+H(x)\tb THP(b)=\Delta(x)\tb b+H(x)\tb THP(b)\\
=(x|x)_H\tb b-H(x)\tb b+TH(x)\tb b+H(x)\tb HP(b)\\
=(x|x)_H\tb b-H(x)\tb b-TTH(x)\tb Tb+H(x)\tb T(b)+H(x)\tb b
=(x|x)_H\tb b.
\end{multline*}
The remaining conditions of $A\tl B$ are trivially satisfied.

\section{Bilinear maps for square groups}
The tensor product of square groups has a universal property similar to that
of abelian groups. To formulate it we need an analog of the notion of a bilinear
map for square groups.

\begin{De}
For square groups $A$, $B$, $C$ a \emph{bilinear map}
$$
\phi:(A,B)\to C
$$
consists of three maps
$$
\phi\ld,\phi\rd:A\e\x B\e\to C\e
$$
and
$$
\phi\ee:A\ee\x B\ee\to C\ee
$$
such that $\phi\ld$ is left linear, $\phi\rd$ is right linear, $\phi\ee$ is
bilinear and moreover one has
\begin{align*}
\phi\ld(P(a),y)&=P\phi\ee(a,\Delta(y)),\\
\phi\rd(x,P(b))&=P\phi\ee(\Delta(x),b),\\
H\phi\ld(x,y)&=\phi\ee((x|x)_H,H(y))+\phi\ee(H(x),\Delta(y)),\\
H\phi\rd(x,y)&=\phi\ee(\Delta(x),H(y))+\phi\ee(H(x),(y|y)_H),\\
P\phi\ee(T(a),T(b))&=-P\phi\ee(a,b),\\
\phi\ld(x,y)-\phi\rd(x,y)&=P\phi\ee(H(x),TH(y)).
\end{align*}

For a bilinear map $\phi$ as above and a morphism $f:C\to C'$ of square groups
the composite bilinear map $f\phi:A\x B\to C'$ is defined by the obvious
equalities
\begin{align*}
(f\phi)\ld(x,y)&=f\e\phi\ld(x,y),\\
(f\phi)\rd(x,y)&=f\e\phi\rd(x,y),\\
(f\phi)\ee(x,y)&=f\ee\phi\ee(x,y).\\
\end{align*}
\end{De}

With this definition we then have

\begin{The}
For square groups $A$, $B$ the rules
\begin{align*}
\un\ld(x,y)&=x\tl y,\\
\un\rd(x,y)&=x\tr y,\\
\un\ee(a,b)&=a\ox b
\end{align*}
define a bilinear map
$$
\un:(A, B)\to A\odot B
$$
which is universal; this means that for any bilinear map $\phi:(A,B)\to C$
there exists a unique morphism $f^\phi:A\odot B\to C$ of square groups with
$\phi=f^\phi\un$.
\end{The}

\begin{proof}
It is straightforward to check that $\un$ indeed is a bilinear map. Moreover
given a bilinear map $\phi$ we define
\begin{align*}
f^\phi\e(x\tl y)&=\phi\ld(x,y),\\
f^\phi\e(x\tr y)&=\phi\rd(x,y),\\
f^\phi\e(a\tb b)&=P\phi\ee(a,b),\\
f^\phi\ee(a\ox b)&=\phi\ee(a,b).
\end{align*}
One checks easily that this indeed defines a morphism of square groups.
Uniqueness is then clear because of the relation $a\tb b=P(a\ox b)$.
\end{proof}

\section{Quadratic functors}

As mentioned in the introduction square groups correspond to quadratic
endofunctors of the category of groups. We now make this correspondence
explicit.

\subsection{Quadratic functors}
Let $\C$ be a category with a zero object and finite coproducts. We choose a
zero object $0\in \C$ and for any objects $X$ and $Y$ we choose a coproduct
$X\vee Y$ in $\C$. We let $i_1:X\to X\vee Y$ and $i_2:Y\to X\vee Y$ be the
structural inclusions. The set $\Hom_\C(X,Y)$ has a distinguished element,
the zero morphism, which is the composite $X\to 0\to Y$ and which is denoted
by $0$ again. For morphisms $f:X\to Z$, $g:Y\to Z$ let $(f,g):X\vee Y\to Z$
be the unique morphism with $(f,g)i_1=f$ and $(f,g)i_2=g$. In particular, we
have the morphisms $r_1:X\vee Y\to X$ and $r_2:X\vee Y\to Y$ given
respectively by $r_1=(\Id_X,0)$ and $r_2=(0,\Id_Y)$. By definition one has
$r_1i_1=\Id_X$, $r_2i_2=\Id_Y$, $r_1i_2=0$ and $r_2i_1=0.$

We consider functors $F:{\C}\to \Gr$ to the category of groups with $F(0)=0$.
For a morphism $f:X\to Y$ in $\C$ let $f_*$ denote the induced morphism
$F(f):F(X)\to F(Y)$. Then the canonical homomorphism
$$
(r_{1*},r_{2*}):F(X\vee Y)\to F(X)\x F(Y)
$$
is always surjective. The kernel of this map is denoted by $F(X|Y)$
and is called the \emph{second cross effect} of $F$. A functor $F$ is called
\emph{linear} provided the second cross effect vanishes, thus $F$ is linear
iff
$$
(r_{1*},r_{2*}):F(X\vee Y)\to F(X)\x F(Y)
$$
is an isomorphism. Moreover, $F$ is called \emph{quadratic} provided $F(X|Y)$
is linear in $X$ and $Y$. Any linear functor has values in the subcategory
$\Ab$ of abelian groups, while values of any quadratic functor lie in the
subcategory $\Nil$ of groups of nilpotence class two \cite{square}.
Moreover, for any quadratic functor $F$ the group $F(X|Y)$ is a central
subgroup of $F(X\vee Y)$ \cite{square}.

\subsection{Pre-square groups and quadratic functors on the category of finite
pointed sets}

Let $\Gamma$ be the category of finite pointed sets. For any $n\ge0$ we
denote by $[n]$ the set $\{0,\cdots,n\}$. We consider $[n]$ as an object of
$\Gamma$ with basepoint $0$. Let $\qgr(\Gamma)$ be the category of all
quadratic functors from $\Gamma$ to $\Gr$. There is an equivalence of
categories
\begin{equation}\label{psq=q}
\PSG\cong\qgr(\Gamma)
\end{equation}
between the category of pre-square groups and the category of quadratic
functors from $\Gamma$ to $\Gr$, which is a particular case of results
obtained in \cite{doldann}. We now discuss functors involved in this
equivalence.

There is a bifunctor
$$
\Gamma\x\PSG\to\Gr
$$
given as follows. Let $M$ be a pre-square group and $S$ a pointed set with basepoint
$*$. Then $S\odot M$ is the group generated by the symbols $s\odot x$ and
$[s,t]\odot a$ with $s,t\in S$, $x\in M\e$, $a\in M\ee$ subject to the
relations
\begin{enumerate}
\item[]
$[s,s]\odot a=s\odot P(a)$
\item[]
$*\odot x=0=[*,s]\odot a $
\item[]
$ [s,t]\odot a= [t,s]\odot T(a)$
\item[]
$[s,t]\odot \{x,y\}=-t\odot x-s\odot y+t \odot x+s\odot y.$
\end{enumerate}
Here $s\odot x$ is linear in $x$ and $[s,t]\odot a$ is central and linear in
$a$. For such $M$ the functor $(-)\odot M:\Gamma\to \Gr$ is quadratic. In
this way one obtains a functor $\PSG\to\qgr(\Gamma)$. A functor $\cro$ in the
opposite direction has the following description. Let $F:\Gamma\to \Gr$ be a
quadratic functor. Then we set
$$\cro(F)\e:=F([1]), \ \ \cro(F)\ee:=F([1]|[1]).$$
Hence one gets a central extension of groups
$$
0\to\cro(F)\ee\to F([2])\to\cro(F)\e\x\cro(F)\e\to 0.
$$
The bracket $\{-,-\}:\cro(F)\e\x\cro(F)\e\to\cro(F)\ee$ is defined by
$$
\{x,y\}:=[F(i_1)(x),F(i_2)(y)]\in F([2])
$$
where $x,y\in F([1])=\cro(F)\e$ and $i_1,i_2:[1]\to [2]$ are pointed maps
given by $i_1(1)=1$ and $i_2(1)=2$ respectively. Moreover, the pointed
involution $[2]\to[2]$ given by $1\mapsto 2$, $2\mapsto 1$ yields an
involution on $\cro(F)\ee$ which is denoted by $T$. We have also the
homomorphism $P:\cro(F)\ee\to\cro(T)\e$ induced by the pointed map $[2]\to
[1]$ given by $1,2\mapsto 1$. One checks that in this way one obtains a
well-defined pre-square group $\cro(F)$. Then the functor
$\cro:\qgr(\Gamma)\to {\PSG}$ is an equivalence of categories whose
quasi-inverse is the functor $M\mapsto (-)\odot M$.

Linear functors correspond
to pre-square groups with $M\ee=0$. Any such functor is isomorphic to a
functor $S\mapsto A\t\bar{\ZZ}(S)$, where $A$ is an abelian group. Here $S$ is a
pointed set with base point $*$, while $\bar{\ZZ}(S)$ is the free abelian group generated by $S$
modulo the relation $*=0$.

\subsection{Square groups and quadratic endofunctors of the category of
groups}\label{kvfunqtorebi}

We restrict ourselves to endofunctors $F:\Gr\to\Gr$ preserving filtered
colimits and reflexive coequalizers. The last condition means that for any
simplicial group $G_*$ the canonical homomorphism $\pi_0(F(G_*))\to
F(\pi_0(G_*))$ is an isomorphism. Such a functor $F$ is completely determined
by the restriction of $F$ to the subcategory of finitely generated free
groups. Let $\lin(\Gr)$ (resp. $\qgr(\Gr)$) be the category of such linear
(resp. quadratic) endofunctors.

A composite of linear endofunctors is linear, therefore $\lin(\Gr)$ is
actually a monoidal category, where the monoidal structure is induced by
composition of endofunctors. Any endofunctor in $\lin(\Gr)$ is isomorphic to
a functor $T$ of the form
$$
T(X)=A\t X\ab
$$
where $A$ is an abelian group. Therefore there is a monoidal equivalence of
monoidal categories
$$
(\lin(\Gr),\circ)\simeq(\Ab,\t).
$$

If $T_i$, $i=1,2$ are quadratic endofunctors of the category of groups, then
the composite $T_1\circ T_2$ in general is not quadratic, but it has a
maximal quadratic quotient which is denoted by $T_2\square T_1$. Then
$\square$ defines a (highly nonsymmetric) monoidal category structure on the
category $\qgr(\Gr)$ \cite{square}.

We now describe an equivalence of categories
$$\Square \cong \qgr(\Gr).$$
For any square group $M$ and any group $G$ one defines the group $G\t M$
\cite{square} by the generators $g\t x$ and $[g,h]\t a$ with $g,h\in G$,
$x\in M\e$ and $a\in M\ee$ subject to the relations
\begin{align*}
(g+h)\t x&=g\t x+h\t x+[g,h]\t H(x),\\
[g,g]\t a&=g\t P(a).
\end{align*}
Here $g\t x$ is linear in $x$ and $[g,h]\t a$ is central and linear in each
variable $g,h$ and $a$. One can check that the functor
$$(-)\t M:\Gr\to \Gr
$$
preserves filtered colimits and reflexive coequalizers.
For any groups $X$ and $Y$ one has by \cite{square} the following
short exact sequence.
$$0\to X\ab\t Y\ab\t M\ee \to (X\vee Y)\t M\to ( X\t M)\x (Y\t M)\to 0$$
Therefore $(-)\t M:\Gr\to \Gr$ is a quadratic functor and hence
this functor is in the category $\qgr(\Gr)$. In this way
 one obtains a functor $\Square\to \qgr(\Gr)$. The functor in the opposite
direction has the following description. If $F:\Gr\to \Gr$ is a quadratic
functor we put
$$
\cro(F)\e=F(\ZZ), \ \ \cro(F)\ee =F(\ZZ\mid\ZZ).
$$
The homomorphism $P$ of the square group $\cro(F)$ is the restriction of the
homomorphism $(\Id,\Id)_*: F(\ZZ *\ZZ)\to F(\ZZ)$. Here $*$ is the coproduct
in the category of groups, so that $\ZZ*\ZZ$ is the free group on two
generators $e_1$ and $e_2$ . The map $H$ is given by
$$
H(x)=\mu_*(x)-p_2(\mu_*x)-p_1(\mu_*x)
$$
Here $\mu:\ZZ\to \ZZ *\ZZ$ is the unique homomorphism which sends $1$ to
$e_1+e_2$, while $p_1$ and $p_2$ are endomorphisms of $\ZZ*\ZZ\to \ZZ*\ZZ$
such that $p_i(e_i)=e_i$, $i=1,2$ and $p_i(e_j)=0$, if $i\not =j$.

For example the quadratic functor corresponding to the square group $A^\t$ is
given by
$$
G\mapsto A\t G\ab\t G\ab.
$$
Here $A\in \Ab$. Similarly, the quadratic functor corresponding to $\Znil[S]$
is given by
$$
G\mapsto\bigvee_{s\in S}G\nilp
$$
where the coproduct is taken in the category $\Nil$.

The main result of \cite{square} shows that
\begin{equation}\label{sq=q}
(\Square, \square)\cong (\qgr(\Gr), \square)
\end{equation}
is a monoidal equivalence of monoidal categories.
Here $\square: \Square\x\Square\to \Square$ is defined
as follows \cite{square}. If $M$ and $N$ are square groups,
then
$$M\square N=(\xymatrix{(M\square N)\e \ar[r]^H &
(M\square N)\ee \ar[r]^P & (M\square N)\e })$$
where $$(M\square N)\e=M\e\t N/ ([x,Pa]\t c\sim 0),
\ \ x\in M\e,a\in M\ee,
c\in N\ee $$
and
$$(M\square N)\ee=((M\ee\t \Cok(P_N))\oplus
(\Cok(P_M)\t \Cok(P_M)\t N\ee))/\sim$$
where one uses the equivalence relation
$$(x|y)_{H_M}\t z\sim \bar{y}\t \bar{x}\t \Delta(z)$$
Moreover, the homomorphism
$P_{M\square N}$ is given by
$$P_{M\square N}(a\t \bar{z})=(Pa)\t z, \ \ P_{M\square N}(\bar{x}\t \bar{y}\t c)=[x,y]\t c$$
while $H_{M\square N}$ is the unique quadratic map satisfying
\begin{align*}
(\bar{x}\t \bar{z}|\bar{x'}\t \bar{z'})_{H_{M\square N}}&=\bar{x'}\t
\bar{x}\t (z|z')_{H_N}\\
H_{M\square N}(x\t z)&=H(x)\t \bar{z}+\bar{x}\t \bar{x}\t H(z),\\
H_{M\square N}([x,y]\t c)&=\bar{x}\t \bar{y}\t c+\bar{y}\t \bar{x}\t T(c).
\end{align*}
The unit object in the monoidal category $(\SG,\square)$ is $\Znil$ from
Section \ref{11}.

\subsection{The functor $\wp:\Square\to \PSG$ in terms of quadratic functors}
Recall that we have the functor
$$
\wp:\Square \to \PSG
$$
defined by
$$
\wp(M)=(M\e,M\ee,T=HP-\Id,(-,-)_H,P),
$$
where $M$ is a square group.

For a pointed set $S$ let $\brk S$ be the free group generated by $S$ modulo
the relation $*=0$, where $*$ is the base point of $S$. Then one has a
natural isomorphism
$$
S\odot \wp(M) \ \cong \ \brk S\tp M.
$$
In other words the following diagram commutes.
$$
\xymatrix{\Gamma\ar[r]^{\langle - \rangle} \ar[rd]_{- \odot \wp(M)}&
\sf{\Gr}\ar[d]^{- \tp M}\\
& \sf{Groups}}
$$
This fact has the following interpretation in the language of functors.
Suppose $F:\Gamma\to\Gr$ and $T:\Gr\to\Gr$ are functors, then the composite
$T\circ F:\Gamma\to\Gr$ is a well-defined functor. Moreover, if $F=\brk-$ is
the free group functor, then $T\circ F$ is quadratic provided $T$ is
quadratic. This follows from the fact that $\brk-$ respects coproducts. Thus
one obtains the functor $\qgr(\Gr)\to \qgr(\Gamma)$ which corresponds to
$\wp$ under the equivalences \eqref{psq=q} and \eqref{sq=q}.

If $F$ and $T$ are quadratic functors then $T\circ F$ in general is not
quadratic (it is of degree $\le4$), but has a maximal quadratic quotient
denoted by $F\square T$. In terms of (pre)square groups this means that there
is a well-defined bifunctor
$$
\square:{\PSG}\x {\SG}\to {\PSG}$$
with the property that for square groups $M$, $N$
$$
\wp(M)\square N=\wp(M\square N).
$$
More precisely, if $M$ is a pre-square group
and $N$ is a square group, then the definition of $M\square N$ mimics
the previous definition in Section \ref{kvfunqtorebi}. The groups $(M\square N)\e$ and $(M\square N)\ee$ and the
homomorphism $P_{M\square N}$ are
defined by the same formul\ae\ , of course now we have to consider the
relation $\{x, y\}_{M}\t z\sim \bar{y}\t \bar{x}\t \Delta(z)$.
The bilinear map $\{-,-\}_{M\square N}$ is defined by
$$
\{\bar x\t \bar z,\bar{x'}\t\bar{z'}\}_{M\square N} =\bar{x'}\t\bar
x\t\cref{H_N}z{z'}
$$
and the involution $T_{M\square N}$ is given by
$$
T_{M\square N}(a\t z)=T_M(a)\t z, \ \ T_{M\square N}(\bar x\t\bar y\t c)=
\bar y\t \bar x\t T_N(c).
$$

\Rem\label{comspectra} It is well known that any functor $\Gamma\to \sf
Spaces$ gives rise to a spectrum (see \cite{segal} and \cite{bfr}). If
$F:\Gr\to \Gr$ is a functor, then it yield a functor $\tilde{F}:\Gamma \to
\sf Spaces$ given by $S\mapsto B(F(\brk S)).$ Here $B(G)$ is the classifying
space of a group $G$. In particular $F$ gives rise to a spectrum ${\sf
sp}(F)$,
%whose $n$-th term is the classifying space of the simplicial group obtained
%by applying $F$ to the Kan loop group of the $(n-1)$-dimensional sphere.
This construction can be applied to $F=(-)\t M$, where $M\in \SG$. According
to \cite{extracta} the homotopy groups of this connective spectrum are
isomorphic to the homology groups of the complex
$$
\cdots
\xto{\Id-T}M\ee\xto{\Id+T}M\ee\xto{\id-T}M\ee\xto{\Id+T}M\ee\xto{\id-T}M\ee\xto
PM\e.
$$
Thus in Corollary \ref{29} we constructed for any square group the
homomorphism $k^M:\pi_0({\sf sp}(F))\to \pi_1({\sf sp}(F))$ for $F=(-)\t M$,
which coincides with the first Postnikov invariant of the spectrum ${\sf
sp}(F)$ (see \cite{extracta}). Similarly, for any commutative quadratic ring
$R$ in Theorem \ref{psioperaciebi} we constructed a nontrivial
cup$_1$-quadratic map $\psi:\pi_0({\sf sp}(F))\to \pi_1({\sf sp}(F))$ for
$F=(-)\t R$. Here $\psi$ is not an invariant of the homotopy type of the
spectrum ${\sf sp}(F)$ but depends on the structure of ${\sf sp}(F)$ given by
the commutative quadratic ring $R$.

\subsection{Square rings}\label{sqrings}

A monoid in the monoidal category $(\SG,\square)$ is termed a \emph{square
ring}. More explicitly a square ring can be defined as follows (see
\cite{BHP}, \cite{square}, \cite{iwa}). A square ring $Q$ is a square group
such that $Q\e$ has additionally a multiplicative monoid structure. The
multiplicative unit of $Q\e$ is denoted by $1$. One requires that this monoid
structure induces a ring structure on the abelian group $\Cok(P_Q)$ through
the canonical projection
$$
Q\e\to \Cok(P_Q), \ \ \ a\mapsto \bar{a}.
$$
Moreover the abelian group $Q\ee$ is a $\Cok(P_Q)\t \Cok(P_Q)\t
(\Cok(P_Q))\op$-module with action denoted by $(\bar t\t\bar s)\cdot
u\cdot\bar{r}\in Q\ee$ for $\bar t, \bar s, \bar{r}\in \Cok(P_Q)$, $u\in
Q\ee$. In addition the following identities must be satisfied, where
$H(2)=H(1+1)$:
\begin{enumerate}
\item[(i)] $(x\mid y )_H=(\bar{y}\t \bar{x})\cdot H(2)$,
\item[(ii)]
$T((\bar{x}\t \bar{y})\cdot a\cdot\bar{z})=(\bar{y}\t \bar{x})\cdot
T(a)\cdot\bar{z}$,
\item[(iii)]
$P(a)x=P(a\cdot x)$,
\item[(iv)]
$xP(a)= P((\bar{x}\t \bar{x})\cdot a)$,
\item[(v)]
$H(xy)=(\bar{x}\t \bar{x})\cdot H(\bar{y})+H(x)\cdot \bar{y}$,
\item[(vi)] $(x+y)z=xz+yz+P((\bar{x}\t \bar{y})\cdot H(z))$,
\item[(vii)] $x(y+z)=xy+xz$.
\end{enumerate}

The category of square rings is denoted by $\sf SR$.
Square rings $Q$ with $Q\ee=0$ are precisely rings. Thus we have a full embedding
$$
\mathsf{Rings}\to\mathsf{SR}
$$
which has a left adjoint given by $Q\mapsto \Cok(P_Q)$. The initial object in
the category of $\sf SR$ is $\Znil$. Observe that in any square ring $Q$ one
has $\Delta(a)=H(2)\cdot\bar a$ for any $a\in Q\e$.

%In terms of functors a square ring structure on a square group $Q$ is
%the same as a monad structure on the functor $(-)\t Q:\Nil\to
%\Nil$. As usual with monads, one can talk about algebras over a monad.
%They are precisely the right modules over $Q$ in the sense of \cite{BHP}.
%For this we recall that a right $Q$-module is a group $G$ together with
%maps $G\x Q\e\to G$, $(g,x)\to ax$ and $G\x G\x Q\ee\to G$,
%$(gmh,a)\mapsto [g,h]_a$ satisfying the following relations.
%\begin{align*}
%m1&=m,\\
%(mx)y&=m(xy),\\
%m(x+y)&=mx+my,\\
%(m+n)x&=mx+nx+[m,n]_{H(a)},\\
%mP(a)&=[m,m]_a,\\
%[m,n]_{Ta}&=[n,m]_a,\\
%[mx,ny]_a&=[m,n]_{(x\otimes y)a},\\
%[[m,n]_a,z]_b&=0.
%\end{align*}
%Moreover $[m,n]_a$ is linear in $m,n$ and $a$.
%Similarly, one can talk about left $Q$-modules. In terms of functors they are
%left $((-)\t Q)$-functors of the form $(-)\t M$ for a square group $M$, or
%equivalently left $Q$-objects in the monoidal category $(\Square, \square)$.
%Let us observe that if $Q=\Znil$, then the category of right
%$\Znil$-modules is nothing but $\Nil$, while the category of left
%$\Znil$-modules is nothing but $\Square$.

\subsection{Relation between the $\square$ and $\tl$ products}

\begin{Pro}\label{squaretl}
For square groups $M$ and $N$ there is a well-defined morphism of square
groups
$$
\sigma=\sigma_{M,N}:M\square N\to M\tl N
$$
given by
\begin{align*}
\sigma\e(x\t z)&=x\tl z,\\
\sigma\e([x,y]\t c)&=(y | x)_H\tb c,\\
\sigma\ee(a\t \bar{z})&=a \t \Delta(z),\\
\sigma\ee(\bar{x}\t \bar{y}\t c)&=(y|x)_H\t c.
\end{align*}
Here $x,y\in M\e$,$z\in N\e$ and $c\in M\ee$. Moreover, $\sigma$ equips the
identity functor
$$
\Id:(\SG,\tl)\to(\SG,\square)
$$
with the structure of a lax monoidal functor.
\end{Pro}

\begin{proof}
First we have to check that $\sigma\e$ is well-defined. In other words
$\sigma\e$ respects all relations of $M\e \t N$ and $\sigma$ vanishes on
$[x,Pa]\t c$. The last assertion is easy to check:
$$
\sigma\e([x,Pa]\t c)=(Pa | x)_H\tb c=0.
$$
We have $\sigma\e((x+y)\t z)=(x+y)\tl z=x\tl z+y\tl z+(y | x)_H\tb H(z)$. On
the other hand $\sigma\e(x\t z+y\t z+[x,y]\t Hz)=x\tl z +y\tl z+(y | x)_H\tb
H(z)$. Thus $\sigma\e$ respects the corresponding
relation of the definition of
the tensor product $M\e\t N$. Similarly
$$
\sigma\e([x,x]\t c)=(x | x)_H\tb c=x\tl P(c),
$$
therefore $\sigma\e$ respects another relation of the definition of the
tensor product $M\e \t N$. Other relations are even easier to check and
therefore we omit them. Next we have to show that $\sigma$ is a morphism
of square groups. We have
$$
\sigma\e P(a\t z)=\sigma\e(Pa\t z)=Pa\tl z=a\tb\Delta(z)=P\sigma\ee(a\t z)
$$
and
$$
\sigma\e P(x\t y\t c)=\sigma\e([x,y]\t c)=(y | x)_H\tb c=P\sigma(x\t y\t c)
$$
which shows compatibility with $P$. To show compatibility with $H$,
first we have to check it for cross-effects
\begin{multline*}
(\sigma\e (x\t z)\mid \sigma\e(y\t z'))_H= (x\tl z\mid y\tl z')\\
=(x\mid y)_H\tl (z\mid z')_H= \sigma\ee (y\t x\t (z\mid z')_H)=\sigma\ee
((x\t z\mid y\t z')_H).
\end{multline*}
Then we have
\begin{multline*}
H(\sigma\e(x\t z))=H(x\tl z)=(x\mid x)_H\t Hz+Hx\t \Delta(z)\\
=\sigma\ee(x\t x\t Hz+Hx\t z)=\sigma\ee(H(x\t z))
\end{multline*}
as well as
\begin{align*}
H\sigma\e([x,y]\t c)&=H((y\mid x)_H\tb c)\\
&=(y\mid x)_H\tb c-T(y\mid x)_H\tb Tc \\
&=\sigma\ee(x\t y\t c+y\t x\t Tc)\\
&=\sigma\ee H([x,y]\t c),
\end{align*}
hence $\sigma$ is compatible with $H$.
%It remains to show that $\sigma$ is coherent, that is
%the following diagram
%$$\xymatrix{& (A\square B)\square C \ar[r]\ar[dl]&A\square (B\square C)\ar[dr]&\\
%(A\square B)\tl C \ar[dr]&&& A\square (B\square C)\ar[dl]\\
%& (A\tl B)\tl C \ar[r] & A\tl (B\tl C)&}
%$$
\end{proof}

Theorem \ref{jibladzestrange} below describes an important case, when the
transformation $\sigma$ is an isomorphism.

\begin{Co} Any quadratic ring $R$ gives rise to a square ring, whose
underlying square group is the same, while the $\Cok(P)\t \Cok(P)\t
\Cok(P)\op$-module structure on $M\ee$ is given by
$$
(\bar{x}\t \bar{y}) a z=(x|y)_Ha\Delta(z).
$$
In this way one obtains a functor
$$
U:\mathsf{QR}\to\mathsf{SR}.
$$
\end{Co}

\section{Abelian groups with cosymmetry as a full monoidal
subcategory}\label{moncosymmetry} In this section we introduce ``abelian
groups with cosymmetry'' which form a subcategory of the category of square
groups. The corresponding quadratic endofunctors of the category of groups of
nilpotence class 2 are coproduct preserving.

\subsection{Symmetric monoidal category $\Cos$}
A pair $(A,\d)$ is called an \emph{abelian group with cosymmetry} if $A$ is
an abelian group and $\d:A\to\Sym^2(A)$ is a map from $A$ to the second
symmetric power of $A$ satisfying
$$\d(a+b)=\d(a)+\d(b)+ab, \ \ a,b\in A.$$

Let $\Cos$ be the category of abelian groups with cosymmetry. We equip this
category with a symmetric monoidal structure. To this end we need the
maps
\begin{align*}
*:A\x\Sym^2(B)&\to\Sym^2(A\t B),\\
*:\Sym^2(A)\x B&\to\Sym^2(A\t B)
\end{align*}
and the homomorphism
$$*:\Sym^2(A)\t\Sym^2(B)\to\Sym^2(A\t B) $$
defined respectively by
\begin{align*}
(a,bb')&\mapsto(a\t b)(a\t b'),\\
(aa',b)&\mapsto(a\t b)(a'\t b),
\end{align*}
and
$$(aa')\t (bb')\mapsto (a\t b)(a'\t b')+(a'\t b)(a\t b').$$
Now we define
$$(A,\d_A)\t (B,\d_B)=(A\t B,\d_{A\t B})$$
where $\d_{A\t B}:A\t B\to\Sym^2(A\t B)$
is given by
$$\d_{A\t B}(a\t b)= a*\d_B(b)+\d_A(a)*b-\d_A(a)*\d_B(b).$$

It is straightforward to check that in this way we really get a symmetric monoidal
category $(\Cos,\t)$, with unit object given by $(\ZZ,\binom -2)$
\subsection{Coproduct preserving endofunctors of the category
$\Nil$} Let $T\in \qgr(\Gr)$ be a functor. By definition $T$ is a quadratic
endofunctor of the category of groups which preserves filtered colimits and
reflexive coequalizers. According to \cite{square} values of $T$ lie in the
category $\Nil$ of groups of nilpotence class two and moreover for any group
$G$ there is an isomorphism
$$T(G)\cong T(G\nilp)$$
This implies that $$\qgr(\Gr)\cong \qgr(\Nil)$$ and we
 identify any object of $\qgr(\Gr)$ with a quadratic endofunctor of
the category $\Nil$
 which preserves filtered colimits and reflexive
coequalizers.

\begin{Pro}\label{copreq} If $F:\Nil\to \Nil$ preserves finite coproducts then
$F$ is quadratic.
\end{Pro}

\begin{proof} Since $F(X\vee Y)\cong F(X)\vee F(Y)$
it follows from the exact sequence \eqref{jaminilsi} that one has the
following exact sequence
$$
0\to F(X)\ab\t F(Y)\ab\to F(X\vee Y)\to F(X)\x F(Y)\to 0.
$$
Thus we have to show that the functor $G$ given by $G(X)=F(X)\ab$ is linear.
Since $(-)\ab:\Nil\to\Ab$ commutes with coproducts, we see that the functor
$G:\Nil\to\Ab$ commutes with finite coproducts, and therefore $G$ is linear,
because finite coproducts in $\Ab$ are finite products as well.
\end{proof}

Let $\qgrs$ be the full subcategory of the category $\qgr(\Nil)$ with finite
coproduct preserving functors as objects. By Proposition \ref{copreq} the
category $\qgrs$ consists of all endofunctors $T:\Nil\to \Nil$ which commute
with all colimits. Therefore we get the following result.

\begin{Co}\label{sq=circ} If $F_i\in \qgrs$, $i=1,2$ then the composite
$F_1\circ F_2\in \qgrs$. Therefore the canonical projection
$$F_2\square F_1\to F_1\circ F_2$$
is an isomorphism.
\end{Co}

It follows that $(\qgrs,\circ)$ is a monoidal category. Below we show
 that it is equivalent as a monoidal category to $(\Cos,\t)$. In particular
 $(\qgrs,\circ)$ is a
symmetric monoidal category.

\begin{Le}\label{crizo} Let $M$ be a square group and let $F=(-)\t M:\Nil\to
\Nil$ be the corresponding functor. Then $F\in \qgrs$ if and only if the
homomorphism of abelian groups
$$(-|-)_H:\Cok( P)\t \Cok(P)\to M\ee$$
is an isomorphism.
\end{Le}

\begin{proof} For the groups $M\e=F(\ZZ)$, $M\ee$ and
$F(\ZZ)\vee F(\ZZ)$ we have the following commutative diagram
with exact rows
$$\xymatrix{0\ar[r]&
M\ee\ar[r] & F(\ZZ\vee \ZZ)\ar[r]& F(\ZZ)\x F(\ZZ)\ar[r]& 0\\
0\ar[r]&
F(\ZZ)\ab\t F(\ZZ)\ab \ar[r]_{[-,-]}\ar[u]^{(-|-)_H}&
 F(\ZZ)\vee F(\ZZ)\ar[r]\ar[u]^{\alpha} &F(\ZZ)\x
 F(\ZZ)\ar[r]\ar[u]_\Id& 0}$$
It follows that the canonical map $\alpha:F(\ZZ)\vee F(\ZZ)\to F(\ZZ\vee
\ZZ)$ is an isomorphism iff $\cref H--:F(\ZZ)\ab\t F(\ZZ)\ab \to M\ee$ is an
isomorphism. Since this map always factors through $\Cok(P)\t \Cok(P)$ it
follows also that the quotient map $F(\ZZ)\ab\t F(\ZZ)\ab\to \Cok(P)\t
\Cok(P)$ is an isomorphism as well. Thus we have proved the ``if'' part.
Assume $\cref H--:\Cok( P)\t \Cok(P)\to M\ee$ is an isomorphism. Then we have
the following exact sequence
$$
0\to M\e\t M\e\t X\ab \t Y\ab\to F(X\vee Y)\to F(X)\x F(Y)\to0.
$$
We have to show that $F(X)\vee F(Y)\to F(X\vee Y)$ is an isomorphism for all
$X,Y\in \Nil$. Since $F$ respects reflexive coequalizers, it suffices to
assume that $X,Y$ are free in $\Nil$. Since $F$ preserves filtered colimits
we can assume that $X$ and $Y$ are finitely generated free. So it suffices to
show that for all $n$ the natural map from the $n$-th copower of $F(\ZZ)$ to
$F(X)$ is an isomorphism, where $X=\ZZ\vee\cdots\vee\ZZ$ ($n$-fold
coproduct). We already proved the statement for $n=2$. Since it is also clear
for $n=0$ or $n=1$ we can proceed by induction. Let $Y$ be the $(n-1)$-st
copower of $\ZZ$. Then we have the following exact sequence
$$
0\to M\e\ab \t M\e\ab\t Y\ab \to F(Y\vee \ZZ)\to F(Y)\x M\e \to 0.
$$
We also have the following exact sequence
$$
0\to F(Y)\ab \t M\e\ab \to F(Y)\vee M\e\to F(Y)\x M\e\to 0.
$$
Since $Y\ab =\ZZ^{n-1}$ and since by the inductive assumption
$F(Y)\ab=(M\e\vee \cdots \vee M\e)\ab=(M\ee \ab)^{n-1}$ the
result follows.
\end{proof}

\begin{Co}\label{cinacos}
The category $\qgrs$ is equivalent to the full subcategory $\Square
_{\Sigma}$ of $\Square$ consisting of square groups $M$ for which
$$
\cref H--:\Cok( P_M)\t \Cok(P_M)\to M\ee
$$
is an isomorphism.
\end{Co}

\begin{proof}
It is enough to notice that $\qgrs\subset \qgr(\Nil)$ and therefore an object
$F\in \qgrs$ is isomorphic to a functor $F=(-)\t M$, with a square group $M$.
The rest follows from Lemma \ref{crizo}.
\end{proof}

Now we define a functor
$$\sJ:\Cos\to \Square _{\Sigma}$$
which is based on a construction from \cite{strange}. Let $\d$ be
a cosymmetry on an abelian group $A$. The square group $\sJ(A,\d)$ is defined
as follows. Consider the pullback diagram in the category of sets:
$$
\xymatrix{
(\sJ(A,\d))\e\ar[r]^-H\ar[d]_p&A\t A\ar[d]^\pi\\
A\ar[r]_-\d&\Sym^2(A)
}
$$
where $\pi(a\t b)=ab$. Then $\sJ(A,\d)\e$ is a group via
$$
(a,x)+(b,y)=(a+b,x+y+a\t b)
$$
where $a,b\in A$ and $x,y\in A\t A$ are such elements that $\pi(x)=\d(a)$ and
$\pi(y)=b$. Let $\sJ(A,\d)\ee$ be $A\t A$, with $P$ given by $P(a\t b)=(0,a\t
b-b\t a)\in \sJ(A,\d)\e$. One easily shows that
$$
\sJ(A,\d)=(\sJ(A,\d)\e\xto H\sJ(A,\d)\ee\xto P\sJ(A,\d)\e)\in \sf SG.
$$

Let $\d$ be a cosymmetry on an abelian group $A$. It follows from the
construction of the square group $M=\sJ(A,d)$ that
$$
\Cok(P_M)\cong A\cong(M\e)\ab
$$
and the commutator map $\Lambda^2(A)\to M\e$ is a monomorphism. Moreover, the map
$$
\cref H--:A\t A\to M\ee
$$
is an isomorphism. Hence we get in fact the functor
$$
\sJ:\Cos\to\Square_\Sigma.
$$
For example, for a set $S$ there is a unique cosymmetry structure on
$A=\ZZ[S]$ such that $\d(s)=0$ for all $s\in S$. In this case one has
$$
\sJ(A,\d)\cong \Znil[S].
$$

The following theorem is a reformulation of a result from \cite{strange}.

\begin{The}\label{jibladzestrange}
\
\begin{itemize}
\item[i)]
The above functor
$$
\sJ:\Cos\to\Square_\Sigma
$$
is an equivalence of categories.
\item[ii)]
A morphism $f:M\to N$ in $\Square_\Sigma$ is an isomorphism iff the induced
homomorphism $\Cok(P_M)\to \Cok(P_N)$ is an isomorphism.
\item[iii)]
Let $M$ and $N$ be square groups. If $M,N\in\Square_\Sigma$, then $M\tl N\in
\Square_\Sigma$ and $M\square N\in\Square_\Sigma$ and the canonical morphism
$\sigma:M\square N\to M\tl N$ is an isomorphism.
\item[iv)]
Let $(A,\d)$ and $(A',\d')$ be abelian groups with cosymmetry. Then there is
a natural isomorphism
$$
\sJ(A,\d)\tl\sJ(A',\d')\cong \sJ((A,\d)\t (A',\d'))
$$
Thus the functor
$$
\sJ:(\Cos,\t ) \to (\SG, \tl)
$$
is symmetric monoidal; moreover it is full.
\end{itemize}
\end{The}

\begin{proof}
By Corollary \ref{cinacos} we know that $\Cos$ is equivalent to the category
of square groups for which $M\ee =\Cok(P) \t \Cok(P)$ and $\{-|-\}:\Cok(P)\t
\Cok(P)\to M\ee$ is the identity map. It follows from the definition of
square groups that
$$
P(a\t b)=a\t b-b\t a
$$
and
$$
H(x+y-x-y)=-b\t a+a\t b.
$$
Thus we have the following commutative diagram with exact rows
$$
\xymatrix{
&\La^2 A\ar[r]^-{\sf c}\ar@{=}[d]&M\e\ar[d]^H\ar[r]&A\ar[r]\ar@{-->}[d]^\delta&0\\
0\ar[r]&\La^2A\ar[r]_-{\sf d}& A\t A\ar[r]_-\pi&\Sym^2(A)\ar[r]&0
}
$$
which yields a cosymmetry $\delta:A\to\Sym^2(A)$. Here $\sf c$ is the
commutator map ${\sf c}(\bar{x}\wedge \bar{y})=x+y-x-y$ and ${\sf d}(a\wedge
b)=a\t b-b\t a$. It follows that ${\sf c}:\La^2 A\to M\e$ is a monomorphism
and $M\cong \sJ(A,\d)$. In this way we have constructed a functor $\Psi$ in
the opposite direction, hence the equivalence of categories in i) is proved.
Explicitly, we have $\Psi(M)=(\Cok(P_M),\delta_M)$, where $\delta_M$ is
induced from the composite
$$
H:M\e\to M\ee\cong \Cok(P_M)\t \Cok(P_M)\to\Sym^2(\Cok(P_M)).
$$
ii) follows from i) because $f$ is an isomorphism iff $\Psi(f)$ is an
isomorphism. iii) First we consider the case of the $\tl$-product. It
suffices to check that $H_{M\tl N}$ yields an isomorphism $\Cok(P_{M\tl N})\t
\Cok(P_{M\tl N})\to (M\tl N)\ee.$ But this is clear, since $\Cok(P_{M\tl N})$
is $\Cok(P_M)\t \Cok(P_N)$. For the $\square$-product this follows from
Corollary \ref{sq=circ} and Corollary \ref{cinacos}. It remains to prove iv),
which claims compatibility with the $\tl$ product. This follows from the
explicit description of this operation.
\end{proof}

\subsection{An open question}
Let $\bf T $ be an algebraic theory, and consider the category ${\sf
End}^{\Sigma}({\bf T})$ of all endofunctors of the category of models of $\bf
T$ preserving all colimits. Then ${\sf End}{(\bf T)}$ is closed with respect
to composition. Under what conditions is the monoidal category $({\sf
End}^{\Sigma}({\bf T}),\circ)$ symmetric? Does this happen if $\bf T$ is the
theory of nilpotent groups of any class $n$? It is obviously so if $n=1$ and
it follows from our results that the same is true for $n=2$. A classical
result of Kan \cite{kan} implies that this also holds for the algebraic
theory of groups.


\begin{thebibliography}{33}

\bibitem{comca} H.-J. Baues. Commutator calculus and
groups of homotopy classes. London Mathematical Society
Lecture Note Series, 50. Cambridge
University Press, Cambridge-New York, 1981. ii+160 pp.

\bibitem{CH} H.-J. Baues. Combinatorial homotopy and
$4$-dimensional complexes. de Gruyter Expositions in Mathematics, 2. Berlin,
1991. xxviii+380 pp.

 \bibitem{qf} H.-J. Baues. Quadratic functors and metastable homotopy.
 J. Pure Appl. Algebra 91 (1994), 49--107.

\bibitem{BHP} H.-J. Baues, M. Hartl and T. Pirashvili. Quadratic
categories and square rings. J. Pure Appl. Algebra 122(1997), 1-40.

\bibitem{iwa} H.-J. Baues and N. Iwase.
Square rings associated to elements in homotopy groups of spheres. Contemt.
Math. 274(2001), 57-78.

\bibitem{Bm} H.-J. Baues and F. Muro. On adding relations to
homotopy groups. (Preprint)

\bibitem{square} H.-J. Baues and T. Pirashvili.
Quadratic endofunctors of the category of groups. Adv. Math. 141 (1999), no. 1,
167--206.

\bibitem{extracta} H.-J. Baues and T. Pirashvili. Spaces associated to
quadratic endofunctors of the category of groups. Extracta Mathematica.
20(2005), 99-136.

\bibitem{shukla} H.-J. Baues and T. Pirashvili. Comparison of Mac Lane,
Schukla and Hochschild cohomologies. (To appear in J. Reine Angew. Math.)

\bibitem{BL} R. Brown and J.-L. Loday. Van Kampen theorems for diagrams of spaces.
With an appendix by M. Zisman. Topology 26 (1987), no. 3, 311--335.

\bibitem{bfr} A. K. Bousfield and E.M. Friedlander. Homotopy
theory of $\Gamma$-spaces, spectra, and bisimplicial sets.
Geometric applications of homotopy theory (Proc. Conf., Evanston, Ill.,
1977), II, pp. 80--130, Lecture Notes in Math., 658, Springer, 1978.

%\bibitem{Day} B. Day. On closed categories of functors. 1970 Reports of the Midwest
%Category Seminar, IV pp. 1--38 Lecture Notes in Mathematics, Vol. 137.

\bibitem{dreck} W. Dreckmann. Distributivgesetze in
der Homotopietheorie. Dissertation, Rheinische Friedrich-Wilhelms-Universität
Bonn, Bonn, 1992. Bonner Mathematische Schriften, 249. Universität Bonn,
Mathematisches Institut, Bonn, 1993. vi+96 pp.
%{igusa} K. Igusa, On the algebraic $K$-theory of $A_{\infty }$-ring
%spaces. Springer Lecture Notes in Math., v. 967 (1982) 146--194.

\bibitem{strange} M. Jibladze. Some strange monoidal
categories. Proc. A. Razmadze Math. Inst. 113 (1995), 83--93.

\bibitem{JP} M. Jibladze and T. Pirashvili. Cohomology
of algebraic theories. J. Algebra 137 (1991), no. 2, 253--296.

\bibitem{niq} M. Jibladze and T. Pirashvili.
Quadratic envelope of the category of class two nilpotent groups. Preprint
MPIM 2005-35.

\bibitem{js} A. Joyal and R. Street. Braided tensor categories.
Adv. Math. 102 (1993), no. 1, 20--78.

\bibitem{kan} D. M. Kan. On monoids and their dual. Bol. Soc. Mat. Mexicana. 3(1958), 52-- 61.

\bibitem{HC} J.-L. Loday. Cyclic homology. Second edition. Grundlehren der Mathematischen
Wissenschaften, 301. Springer-Verlag, Berlin, 1998. xx+513 pp.

\bibitem{working} S. MacLane. Categories for the working mathematician. Graduate
Texts in Mathematics, Vol. 5. Springer-Verlag, New York-Berlin, 1971. ix+262 pp.

\bibitem{doldann} T. Pirashvili. Dold-Kan
type theorem for $\Gamma$-groups. Math. Ann. 318 (2000), no. 2, 277--298.

\bibitem{PW} T. Pirashvili and F. Waldhausen. Mac Lane homology
and topological Hochschild homology. J. Pure Appl. Algebra 82 (1992), 81--98.

\bibitem{qss} D. G. Quillen. Spectral sequences of a double semi-simplicial group.
 Topology 5 1966 155--15.

\bibitem{qu} D. G. Quillen. Homotopical algebra. Lecture
Notes in Mathematics, No. 43 Springer-Verlag, Berlin-New York 1967 iv+156 pp.

\bibitem{segal} G. Segal. Categories and cohomology theories.
Topology 13 (1974), 293--312.

\end{thebibliography}
\end{document}